\documentclass[a4paper, 11pt]{article}

\usepackage{amsmath}
\usepackage{amssymb}
\usepackage{cite}
\usepackage{makeidx}
\usepackage{graphicx}
\usepackage{multirow}
\usepackage{booktabs}
\usepackage{appendix}
\usepackage{framed}
\usepackage[normalem]{ulem}
\usepackage{amsthm}
\usepackage{subfig}
\usepackage{diagbox}
\usepackage{algorithm}
\usepackage{color}
\usepackage{bigstrut}
\usepackage{algpseudocode}
\usepackage{amsfonts}
\usepackage{grffile}
\usepackage[colorlinks,linkcolor=blue,anchorcolor=blue,citecolor=red]{hyperref}
\usepackage{enumerate}
\numberwithin{equation}{section}
\newcommand{\UP[1]}{^{(#1)}}
\newcommand\pd[2]{\frac{\partial {#1}}{\partial {#2}}}
\newcommand\pa{\partial}
\newcommand\rd{\mathrm{d}}
\newcommand\bA{\boldsymbol{A}}
\newcommand\bu{\boldsymbol{u}}
\newcommand\bg{\boldsymbol{g}}
\newcommand\bh{\boldsymbol{h}}

\newcommand\bi{\boldsymbol{i}}
\newcommand\bj{\boldsymbol{j}}
\newcommand\bv{\boldsymbol{v}}
\newcommand\bw{\boldsymbol{w}}
\newcommand\bx{\boldsymbol{x}}
\newcommand{\bsigma}{\boldsymbol{\sigma}}
\newcommand{\bq}{\boldsymbol{q}}
\newcommand\bdf{\boldsymbol{f}}
\newcommand\bF{\boldsymbol{F}}
\newcommand\bT{\overline{\theta}}
\newcommand\oT{\overline{\theta}}
\newcommand\ou{\overline{\bu}}
\newcommand\ot{\overline{\te}}
\newcommand\bQ{\boldsymbol{Q}}
\newcommand{\bz}{\boldsymbol{0}}

\newcommand\bbR{\mathbb{R}}
\newcommand\bbN{\mathbb{N}}
\newcommand\mO{\mathcal{O}}
\newcommand\mQ{\mathcal{Q}}
\newcommand\mR{\mathbb{R}}

\newcommand\mM{\mathcal{M}}
\newcommand\mH{\mathcal{H}}

\newcommand\mZ{\mathbb{Z}}

\newcommand\sv{\bv_{\ast}}
\newcommand\al{\alpha}
\newcommand\be{\beta}

\newcommand{\eps}{\varepsilon}
\newcommand\si{\sigma}
\newcommand\te{\theta}

\newcommand\sk{\substack}

\newcommand\aut{_{\al}^{[\ou,\ot]}}

\newcommand\but{_{\be}^{\ou,\ot}}

\newcommand\Ha{H_{\al}}
\newcommand{\mC}{\mathcal{C}}

\newcommand\mHk{\mH_{\ka}}
\newcommand\mHl{\mH_{\la}}
\newcommand\Hl{H_{\la}}

\newcommand\Hk{H_{\ka}}

\newcommand{\mS}{\mathcal{S}}
\newcommand{\mV}{\mathcal{V}}
\newcommand\ga{\gamma}
\newcommand{\om}{\omega}

\newcommand\de{\delta}
\newcommand\la{\lambda}

\newcommand\ka{\kappa}

\newcommand{\Kn}{{\rm Kn}}
\newcommand\alk{_{\alpha,\lambda,\kappa}}
\newcommand\pq{\preceq}
\newcommand\ca[1]{a_{l'_{#1}k'_{#1}}^{l_{#1}k_{#1}}}

\newcommand\Aalk{A_{\alpha, \lambda, \kappa}}
\newtheorem*{definition}{Definition}
\newtheorem{property}{Property}
\newtheorem{theorem}{Theorem}
\newtheorem{lemma}{Lemma}
\newtheorem{proposition}[theorem]{Proposition}
\newtheorem{corollary}{Corollary}
\newtheorem{remark}{Remark}

\graphicspath{{image/}}

\title{Hermite spectral method for the inelastic Boltzmann equation}
\author{
Ruo Li\thanks{CAPT, LMAM \& School of Mathematical Sciences,
    Peking University, Beijing, China, email: {\tt
      rli@math.pku.edu.cn}.},~~Yixiao Lu\thanks{HEDPS, Center for Applied Physics and Technology \& School of Mathematical
    Sciences, Peking University, Beijing, China, 100871, email: {\tt
      luyixiao@pku.edu.cn}.},~~Yanli Wang\thanks{Beijing Computational
    Science Research Center, email: {\tt ylwang@csrc.ac.cn}.}}
\begin{document}
\maketitle
\begin{abstract}
We propose a Hermite spectral method for the inelastic Boltzmann equation, which makes two-dimensional periodic problem computation affordable by the hardware nowadays. This new approach involves utilizing a Hermite expansion, whereby the expansion coefficients for the VHS model are simplified into a series of summations that can be precisely derived. Additionally, a new collision model is built with a combination of the quadratic collision operator and a simplified collision operator, which helps us to balance the computational cost and the accuracy. Various numerical experiments, including spatially two-dimensional simulations, demonstrate the accuracy and efficiency of this numerical scheme.
\vspace*{4mm}
    \newline
    \noindent \textbf{Keywords:} granular gas flow; 
    inelastic Boltzmann equation; Hermite spectral method
\end{abstract}

\section{Introduction}
\label{sec:intro}
Recently, there has been increasing interest in studying granular materials such as sand, grains, and snow. Unlike molecular particles, which are typically modeled with elastic collision, granular gases exhibit distinct behavior due to the dissipation of energy during collisions. Therefore, most theories for the elastically colliding spheres are insufficient to describe the granular gases. The Boltzmann equation, an important model in elastic theory, can also be extended to effectively capture the behavior of granular gases. Furthermore, the inelastic Boltzmann equation has found applications in modeling social and biological systems \cite{Pareschi2014}.


Due to energy loss, the inelastic collision operator is fundamentally different from the elastic operator. Both analytical and numerical theories in this field are still at an early stage, and we refer readers to the recent reviews \cite{Villani2006, Carrillo2021} for some related results and open questions. Numerically, some methods have been proposed to solve the inelastic Boltzmann equation. The Direct Simulation Monte Carlo (DSMC) method \cite{Bird}, initially developed for the elastic Boltzmann equation, has recently been extended to the inelastic case \cite{Gamba2005, Astillero2005}. It can efficiently simulate the highly rarefied situations but does not work well in low-speed and unsteady flows. In recent years, deterministic methods have made significant progress in kinetic theory. For example, the Fourier spectral method \cite{Hu2016, Mouhot} has been applied to simulate the Boltzmann equation, and subsequently extended to the inelastic case \cite{Filbet2005, Hu2019, Wu2015}. Additionally, the Petrov-Galerkin spectral method has been proposed for the inelastic Boltzmann equation \cite{Hu2020}, and a unified gas-kinetic scheme has been adopted to handle the inelastic collision of granular gases \cite{Liu2019}.

In the study of inelastic gas flows, particular attention is often devoted to the behavior of macroscopic variables, especially temperature. Therefore, we focus on the Hermite spectral method, which allows us to express important macroscopic variables such as density and temperature, using expansion coefficients up to the first few orders. The history of the Hermite spectral method can be traced back to Grad's work \cite{Grad} in 1949, which is known as the moment method. It relies on the concept of using the steady state Maxwellian as the weight function. The distribution function is then expanded using orthogonal polynomials with this weight function, which in this case are Hermite polynomials. In the past few years, remarkable progress has been achieved in applying the Hermite spectral method to solve the Boltzmann equation. An algorithm to approximate the general quadratic Boltzmann collision operator was first derived in \cite{Approximation2019}. Subsequently, the method was verified with the success in the simulation of rarefied gas flow \cite{ZhichengHu2019, multi2022}. Furthermore, it has been modified and extended to address the Vlasov-type equations \cite{FPL2020, Blaustein2023, Bessemoulin2022}.

In this paper, we develop a numerical algorithm based on the Hermite spectral method to solve the inelastic Boltzmann equation. Although the Maxwellian is no longer the steady state of inelastic collisions, the Hermite spectral method retains the advantage of providing straightforward expressions for macroscopic variables. Moreover, the steady state for the inelastic collision has the form of a Dirac-distribution \cite{Filbet2005}, which can also be approximated by a Gaussian function by appropriately choosing the scaling factor. Using this Gaussian function as the weight function, it is expected that the distribution function can be approximated using the corresponding orthogonal polynomials. In the simulations, the scaling factor is chosen as the local macroscopic temperature when approximating the inelastic collision operator. We first derive the algorithm of the inelastic quadratic collision term within the framework of the Hermite spectral method, significantly reducing the computational complexity of calculating the expansion coefficients. For the VHS model, these coefficients can even be exactly obtained through several summations. Next, to balance the accuracy and computational cost, a new collision model is proposed by combining the quadratic collision model with a simplified inelastic collision model modified from previous work \cite{Filbet2013, Astillero2005}. Following the approach in \cite{ZhichengHu2019, multi2022}, we utilize the Strang splitting method to separate the convection and collision parts. The finite volume method is employed to solve the convection term similarly to \cite{ZhichengHu2019}. The collision term can be efficiently computed using the new collision model, which greatly reduces computational costs while maintaining reliable numerical accuracy.

In the numerical experiments, two important spatially homogeneous experiments are first implemented in granular gas flow, including the heating source problem \cite{Noije1998} and Haff's cooling law \cite{Haff1983}. Then tests are conducted on one-dimensional benchmark problems, including Couette flow and Fourier heat transfer. Finally, a two-dimensional periodic diffusion is simulated to further validate the accuracy and efficiency of the method. The numerical solution shows excellent agreement with the reference solution obtained from the direct simulation Monte Carlo (DSMC) method.

The rest of this paper is organized as follows: in Sec. \ref{sec:pre}, we introduce the inelastic collision operator and the general framework of the Hermite spectral method. Sec. \ref{sec:method} describes the algorithm for discretizing the collision term, along with special simplifications for the VHS model. The complete numerical scheme is given in Sec. \ref{sec:numerical}, followed by presentation of the numerical experiments in Sec. \ref{sec:experiment}. The paper ends with some concluding remarks in Sec. \ref{sec:conclusion} and several supplementary contents in the Appendix.  
\section{Inelastic Boltzmann equation and Hermite spectral method}
\label{sec:pre}
In this section, we will first provide a brief review of the Boltzmann equation and the inelastic collision model, and then introduce the general framework for solving the Boltzmann equation using the Hermite spectral method. 
\subsection{Inelastic Boltzmann equation}
\label{sec:Boltzmann}
The behavior of inelastic gas flow can be described by the general form of the Boltzmann equation as follows \cite{Brilliantov}:
\begin{equation}
    \label{eq:Boltzmann}
     \pd{f}{t}+\bv\cdot\nabla_{\bx} f=\frac{1}{\Kn}\mQ[f,f](\bv),
\end{equation}
where $f(t, \bx, \bv)$ is the distribution function, depending on time $t\in\bbR^+$, physical space $\bx\in\bbR^3$ and particle velocity $\bv\in\bbR^3$. $\Kn$ is a non-dimensionalization parameter that reflects the gas properties and the reference length. The nonlinear quadratic collision operator $\mQ$ models the effects of inelastic collisions. 

When particles with velocities $(\bv, \sv)$ collide, the post-collision velocity pair $(\bv',\sv')$ can be expressed using the $\sigma$-representation \cite{Brilliantov, Carrillo2021}:
\begin{equation}
    \label{eq:post-velo}
    \left\{
    \begin{array}{l}
        \bv'=\frac{\bv+\sv}{2}+\frac{1-e}4
        (\bv-\sv)+\frac{1+e}{4}|\bv-\sv|\sigma,\\
        \sv'=\frac{\bv+\sv}{2}-\frac{1-e}4
        (\bv-\sv)-\frac{1+e}{4}|\bv-\sv|\sigma,\\
    \end{array}
    \right.
\end{equation}
where $\sigma$ is a unit vector in $S^2$, and $e\in[0,1]$ represents the restitution coefficient. 
During collisions, the conservation of momentum can be derived as
\begin{equation}
    \label{eq:conserv}
    \bv+\sv=\bv'+\sv', 
\end{equation}
and the dissipation of energy is given by
\begin{equation}
    \label{eq:energy_loss}
    |\bv|^2+|\sv|^2-(|\bv'|^{2}+|\sv'|^{2})=\frac{1-e^2}{4}|\bg|(|\bg|-\bg\cdot\sigma),
\end{equation}
where $\bg = \bv - \sv$ represents their relative velocity. The specific weak form of $\mQ$ in the $\sigma$-representation can be expressed as \cite{Filbet2005, Hu2019}
\begin{equation}
    \label{eq:weak_1}
    \begin{split}
        &\int_{\mR^3}\mQ[f,f](\bv)\phi(\bv)\rd \bv\\
        &\qquad =\frac12\int_{\mR^3}\int_{\mR^3}\int_{S^2}B(|\bg|,\sigma)f(\bv)f(\sv)
        \big[\phi(\bv')+\phi(\bv'_{\ast})-\phi(\bv)-\phi(\bv_{\ast})\big]\rd\sigma\rd\bv\rd \sv,
    \end{split}
\end{equation} 
or
\begin{equation}
    \label{eq:weak_2}
    \begin{split}
        \int_{\mR^3}\mQ[f,f](\bv)\phi(\bv)\rd \bv
         = \int_{\mR^3}\int_{\mR^3}\int_{S^2}B(|\bg|,\sigma)f(\bv)f(\sv)
         \big[\phi(\bv')-\phi(\bv)\big]\rd\sigma\rd\bv\rd \sv,
    \end{split}
\end{equation} 
where $\phi(\bv)$ is a suitable test function such that $\int_{\mR^3}\mQ[f,f](\bv)\phi(\bv)\rd \bv$ is integrable. 

\begin{remark}
    The strong form of the inelastic collision operator $\mQ$ can be derived using the reflection map in the $\omega$-representation, where $\omega$ represents the impact direction. We refer the readers to \cite{Carrillo2021} for more details. 

    For the coefficient $e$, although it always depends on the relative velocity $|\bg|$ in realistic scenarios, we only consider the constant case in the simulations. More discussions on $e$ can be referred to \cite{Carrillo2021} and the references therein.
\end{remark}

In \eqref{eq:weak_1} and \eqref{eq:weak_2}, $B$ is the collision kernel which depends on the type of interactions. The most commonly used form for the inelastic case is the variable hard sphere (VHS) model \cite{Bird}:
\begin{equation}
    \label{eq:VHS}
    B(|\bg|,\sigma)=C_{\varpi}|\bg|^{2(1-\varpi)},
\end{equation}
where $C_{\varpi}>0$ and $0.5\leqslant\varpi\leqslant 1$ are constants. Especially, the Maxwell molecules and the hard sphere model (HS) correspond to $\varpi = 1$ and $\varpi=0.5$, respectively.

Finally, the steady-state solution satisfying $\mQ(\mM, \mM) = 0$ takes the form of a locally Dirac-distribution \cite{Filbet2005}, given by
\begin{equation}
\label{eq:Maxwellian}
    \mM = \delta_{\rho, \bu}(\bv) = \rho \delta (\bv - \bu),
\end{equation}
where $\rho$ and $\bu$ are the density and macroscopic velocity 
\begin{equation}
     \label{eq:macro_1}
     \begin{aligned}
     &\rho(t,\bx) = \int_{\mR^3}f(t,\bx,\bv)\rd\bv, \qquad 
    \bu(t,\bx)=\frac{1}{\rho}\int_{\mR^3}\bv f(t,\bx,\bv)\rd\bv.
     \end{aligned}
\end{equation}


Due to the complexity of the original quadratic collision operator, simplified collision models have been proposed to approximate it. In the elastic case, simplified models such as the BGK model \cite{BGK} have been developed to simplify the original quadratic operator. However, these models do not perform well in the inelastic case where total energy is not conserved during collisions. In \cite[(3.12)]{Filbet2013}, a simplified collision operator with energy dissipation is constructed with
\begin{equation}
    \label{eq:linear}
    \mQ_{\rm sim}(\bv)= \nu_1
    \left[f_{G}(\bv)-f(\bv)\right]+
    \nu_2 \nabla_{\bv}\cdot\left[(\bv-\bu)f(\bv)\right],
\end{equation}
where $\nu_1$ and $\nu_2$ are problem-dependent parameters. In \cite{Filbet2013}, $\nu_1 = \tau(P)$, where $\tau(\cdot)$ is a given function depending on the kinetic pressure $P:=\rho\theta$ and $\nu_2 = C \rho \theta^{1/2} $ with $C$ being a given constant. A similar model is also proposed in \cite[(2.20)]{Astillero2005}, which includes the same energy loss term but with a different $\nu_2$. 
The first term in \eqref{eq:linear} corresponds to the ellipsoidal statistical BGK (ES-BGK) operator, which has the form \cite{Holway}
\begin{equation}
    \label{eq:ESBGK}
    \begin{split}
        f_{G}&=\frac{\rho}{\sqrt{|2\pi \Lambda|}}\exp\left(-\frac12 (\bv-\bu)^T\Lambda (\bv-\bu)\right), \\
        \Lambda&=(\lambda_{ij})\in \bbR^{3\times 3}, \quad \lambda_{ij}=\frac{1}{\Pr}\theta\delta_{ij}+\left(1-\frac{1}{\Pr}\right)\frac{\sigma_{ij}+\delta_{ij}\rho\theta}{\rho},
    \end{split}
\end{equation}
where $\delta_{ij}$ is the Kronecker delta, and $\Pr$ is the Prandtl number, which takes the value of $\frac23$ for monatomic gases. The temperature and stress tensor are denoted by $\theta$ and $\bsigma$, respectively, and they are related to the distribution function as follows:
\begin{equation}
     \label{eq:macro}
     \begin{aligned}
    &3\theta(t,\bx)\rho(t, \bx)= \int_{\mR^3}|\bv-\bu|^2f(t,\bx,\bv)
    \rd\bv, \\
     &\bsigma(t,\bx)=\int_{\mR^3}\left[(\bv-\bu)\otimes(\bv-\bu)-\frac13|\bv-\bu|^2I\right] f(t,\bx,\bv)\rd\bv. 
     \end{aligned}
\end{equation}
The heat flux $\bq$ can be derived from the distribution function as well:
\begin{equation}
    \label{eq:q}
    \bq(t,\bx) = \frac12 \int_{\mR^3}(\bv-\bu)|\bv-\bu|^2 f(t,\bx,\bv)\rd\bv.
\end{equation}

So far, we have introduced the inelastic Boltzmann equation and discussed its properties. Several numerical methods have been developed to tackle this equation, such as the DSMC method \cite{Gamba2005, Astillero2005}, Fourier spectral method \cite{Filbet2005, Hu2019, Wu2015}, and Petrov-Galerkin spectral method \cite{Hu2020}. In this work, a numerical scheme will be developed based on the Hermite spectral method, which offers higher efficiency in capturing the evolution of temperature.

\subsection{Hermite spectral method}
\label{sec:Hermite}
This section presents the general framework for solving the inelastic Boltzmann equation using the Hermite spectral method. The first step is to choose a weight function and then utilize orthogonal polynomials as the basis functions. Precisely, with an expansion center $[\ou, \oT] \in \bbR^3 \times \bbR^+$, the weight function takes the form
\begin{equation}
    \label{eq:maxwellian_expan}
      \omega_{[\ou, \oT]}(\bv) =  \frac{1}{(2\pi \oT)^{\frac32}}
    \exp\left(-\frac{|\bv-\ou|^2}{2 \oT}\right),
\end{equation}
and the corresponding Hermite polynomials are defined as 
\begin{definition}[Hermite Polynomials]
\label{def:Her}
For $\alpha=(\al_1,\al_2,\al_3) \in \bbN^3$, the three-dimensional 
Hermite polynomial $H_{\alpha}^{[\ou,\oT]}(\bv)$ is defined as follows:
\begin{equation}
    \label{eq:Hermite}
    H_{\alpha}^{[\ou,\oT]}(\bv)=\frac{(-1)^{|\alpha|}\oT^{\frac{|\alpha|}{2}}}
    {\omega_{\ou,\oT}(\bv)}
    \dfrac{\partial^{|\alpha|}}{\partial \bv^{\alpha}}
    \omega_{\ou,\oT}(\bv), 
\end{equation}
where $|\alpha|=\al_1+\al_2+\al_3$ and $\partial \bv^{\alpha} = \pa v_1^{\al_1}\pa v_2^{\al_2}\pa v_3^{\al_3}$. The Hermite polynomials possess several useful properties when approximating the complex collision term, which are listed in Appendix \ref{app:Her}.
\end{definition}

Following a similar routine as in \cite{Approximation2019}, one can approximate the distribution function $f(t, \bx, \bv)$ as 
\begin{equation}
    \label{eq:Her-expan}
    f(t,\bx,\bv) \approx \sum_{|\alpha| \leqslant M}f^{[\ou,\oT]}_{\alpha}(t,\bx)
		\mH\aut(\bv),
\end{equation}
where $\mH\aut(\bv)=H\aut(\bv)\omega_{\ou, \oT}(\bv)$ are the basis functions, and $M \in \mathbb{N}$ is the expansion order. The expansion coefficients $f^{[\ou,\oT]}_{\alpha}(t,\bx)$ can be obtained using the orthogonality of the basis functions \eqref{eq:orth}:
\begin{equation}
    \label{eq:falpha}
    f^{[\ou,\oT]}_{\alpha}(t,\bx)=\frac{1}{\alpha!}\int_{\mR^3}f(t,\bx,\bv) H\aut(\bv)
    \rd\bv.
\end{equation}
With the expansion \eqref{eq:Her-expan}, the macroscopic variables in \eqref{eq:macro}, \eqref{eq:q} can be expressed in terms of $f^{[\ou,\oT]}_{\alpha}$ as
\begin{equation}
    \label{eq:macro_f}
    \begin{split}
    &\rho = f^{[\ou,\oT]}_{0},\quad u_k= \overline{u}_k+\frac{\sqrt{\oT}}{\rho}f^{[\ou,\oT]}_{e_k}, 
    \quad \theta=\frac{2\oT}{3\rho}\sum_{k=1}^3f^{[\ou,\oT]}_{2e_k}+\oT-\frac{1}{3\rho}|\bu-\ou|^2, \\
    &\sigma_{kl}=(1+\delta_{kl})\oT f^{[\ou,\oT]}_{e_i+e_j}+\delta_{kl}\rho\left(\oT - \theta\right)-
    \rho\left(\overline{u}_k-u_k\right)\left(\overline{u}_l-u_l\right), \\
    &q_k=2\oT^{\frac32}f^{[\ou,\oT]}_{3e_k}+(\overline{u}_k-u_k)\oT f^{[\ou,\oT]}_{2e_k}+|\ou-\bu|^2\sqrt{\oT} f^{[\ou,\oT]}_{e_k} +\\ 
    & \qquad \qquad\sum_{l=1}^3\left[\oT^{\frac32}f^{[\ou,\oT]}_{2e_l+e_k}+\left(\overline{u}_l-u_l\right)\oT f^{[\ou,\oT]}_{e_l+e_k}+  
     \left(\overline{u}_k-u_k\right)\oT f^{[\ou,\oT]}_{2e_l}\right], \qquad k, l = 1, 2, 3,
    \end{split}
\end{equation}
where $e_1=(1,0,0), e_2=(0,1,0), e_3=(0,0,1)$ represent the unit vectors. 
Therefore, the macroscopic quantities can be easily obtained under the framework of the Hermite spectral method. This allows us to accurately govern the evolution of important macroscopic variables even with a small expansion order $M$.

Assume the collision term is also expanded and approximated using the same basis functions
\begin{equation}
    \label{eq:expan_coll}
    \mQ[f,f](\bv)\approx\sum_{|\alpha|\leqslant M}Q^{[\ou,\oT]}_{\alpha}(t, \bx)\mH\aut(\bv),\qquad  Q^{[\ou,\oT]}_{\alpha}(t, \bx)=\frac{1}{\alpha!}\int_{\mR^3}\mQ[f,f](\bv)
    H\aut(\bv)\rd \bv.
\end{equation}
By substituting the expansion \eqref{eq:Her-expan}, \eqref{eq:expan_coll} into the Boltzmann equation \eqref{eq:Boltzmann}, and matching the coefficients on both sides, one can derive the moment equations as
\begin{equation}
    \label{eq:moment}
    \pd{}{t}f^{[\ou,\oT]}_{\alpha}+\sum_{d=1}^3\pd{}{x_d}\left((\alpha_d+1)\sqrt{\ot}f^{[\ou,\oT]}_{\alpha+e_d}+\overline{u}_df^{[\ou,\oT]}_{\alpha}+\sqrt{\ot}f^{[\ou,\oT]}_{\alpha-e_d}\right)=Q^{[\ou,\oT]}_{\alpha},\quad |\alpha|\leqslant M,
\end{equation}
where the recurrence relationship \eqref{eq:recur} is utilized to handle the convection term. In \eqref{eq:moment}, $f^{[\ou,\oT]}_{\alpha}$ is regarded as $0$ if $\alpha$ contains any negative index or $|\alpha|>M$.

Until now, we have derived the moment equations for the Boltzmann equation. In fact, the evolution of macroscopic variables can be precisely governed by these moment equations due to the relationships \eqref{eq:macro_f}. The main challenge of solving these equations lies in approximating $Q^{[\ou,\oT]}_{\alpha}(t, \bx)$ in \eqref{eq:expan_coll}, which will be discussed in detail in the following section.

\begin{remark}
    It is worth noting that the computational cost can be greatly reduced with a properly chosen expansion center. Furthermore, it is possible to choose different expansion centers in different computational steps \cite{FPL2020, ZhichengHu2019}. 
    
    In the case of the classical Boltzmann equation, the expansion center is chosen based on local macroscopic velocity and temperature to approximate the quadratic collision term \cite{Approximation2019}. This involves utilizing the local Maxwellian
\begin{equation}
    \label{eq:maxwellian}
   \mM_{[\bu, \theta]}(\bv) =  \rho \omega_{[\bu, \theta]}(\bv) =  \frac{\rho}{(2\pi \theta)^{\frac32}}
    \exp\left(-\frac{|\bv-\bu|^2}{2 \theta}\right),
\end{equation}
which represents the steady-state solution, as the weight function to generate the basis polynomials. On the other hand, during the convection step, the expansion center is often chosen as a rough average of the entire domain to ensure numerical accuracy and stability \cite{ZhichengHu2019}.
\end{remark}

\section{Approximation of the collision terms}
\label{sec:method}
In this section, we will focus on the approximation of the collision term within the framework of the Hermite spectral method.
The discretization of the quadratic term will be presented in Sec. \ref{sec:coll}, while the simplification of the VHS model will be covered in Sec. \ref{sec:VHS}.





\subsection{Series expansion of general collision terms}
\label{sec:coll}
Let us first discuss the algorithm to calculate the expansion coefficients $Q^{[\ou,\oT]}_{\alpha}(t, \bx)$ of the quadratic collision term in \eqref{eq:expan_coll}.
With the weak form \eqref{eq:weak_2}, one can simplify \eqref{eq:expan_coll} as
\begin{equation}
    \label{eq:simp1_Q}
    Q^{[\ou,\oT]}_{\alpha}(t, \bx)
    =\int_{\mR^3}\int_{\mR^3}\int_{S^2}
    B(|\bg|,\sigma)f(\bv)f(\sv)\big[H\aut(\bv')-H\aut(\bv)
    \big]\rd\sigma\rd\bv\rd \sv.
\end{equation}
Substituting the expansion of the distribution function \eqref{eq:Her-expan} into \eqref{eq:simp1_Q}, one can derive that 
\begin{equation}
    \label{eq:simp2_Q}
     Q^{[\ou,\oT]}_{\alpha}(t, \bx) = \sum_{|\lambda|\leqslant M}\sum_{|\kappa| \leqslant M}A^{[\ou,\oT]}\alk f^{[\ou,\oT]}_{\lambda}f^{[\ou,\oT]}_{\kappa},
\end{equation}
where
\begin{equation}
    \label{eq:Aalk}
	\begin{split}
	A^{[\ou,\oT]}\alk=&\frac{1}{\alpha!}\int_{\mR^3}\int_{\mR^3}\int_{S^2}B
	(|\bg|,\sigma)
	\mHl^{[\ou,\oT]}(\bv)\mHk^{[\ou,\oT]}(\sv)[H\aut(\bv')-H\aut(\bv)]
    \rd\sigma \rd\bv \rd\sv.
	\end{split}
\end{equation}
Before introducing the algorithm to calculate $A^{[\ou,\oT]}\alk$, we propose Proposition \ref{thm:change_center} below 
\begin{proposition}
\label{thm:change_center}
     In the calculation of $\Aalk^{[\ou, \oT]}$, the coefficients satisfy the relationship
     \begin{equation}
        \label{eq:VHS_A}
        \Aalk^{[\ou,\oT]}=\oT^{1-\varpi}\Aalk^{[\bz, 1]}.
    \end{equation}
    For more details, readers can refer to \cite[Sec. 3.1]{multi2022}. Thus, it is sufficient to compute and store $\Aalk^{[\bz, 1]}$ with expansion center $[\bz,1]$.
\end{proposition}

Consequently, we assume the expansion center to be $[\ou, \ot]=[\bz, 1]$ and omit the superscripts $[\ou, \ot]$ as follows: 
{\small
\begin{equation}
    \label{eq:short}
    f_{\alpha}(t, \bx) = f^{[\bz, 1]}_{\alpha}(t,\bx), \qquad  Q_{\alpha}(t,\bx) = Q^{[\bz, 1]}_{\alpha}(t,\bx), \qquad H_{\alpha}(\bv) = H^{[\bz, 1]}_{\alpha}(\bv), \qquad \omega(\bv) = \omega_{[\bz, 1]}(\bv). 
\end{equation}
}

With the properties of Hermite polynomials, the calculation of $A\alk$ can be greatly simplified, and the result is listed in Theorem \ref{thm:step1}. 
\begin{theorem} \label{thm:step1}
    The coefficients $A\alk$ have the following form:
\begin{equation}
    \label{eq:theorem_A}
    A\alk=\frac{1}{2^{\frac{|\alpha|+3}{2}}}
    \sum_{\lambda'\pq\alpha,\lambda'\pq\kappa+\lambda}\frac{1}{\bj!}
    \ca1\ca2\ca3\ga_{\kappa'}^{\bj},
\end{equation}
where $\kappa'=\kappa+\lambda-\lambda'$, $\bj=\alpha-\lambda'$, and the symbol `$\pq$' means
$$(i_1,i_2,i_3) \pq (j_1,j_2,j_3) \Leftrightarrow 
i_s\leqslant j_s \quad\text{for}\quad s=1,2,3.$$
The coefficient $\ca{d}$ is computed with
\begin{equation}
    \label{eq:ca}
    \ca{d}=2^{-\frac{l'_d+k'_d}{2}}
    \sum_{s\in \mZ}C_{l_d}^sC_{k_d}^{l'_d-s}
    (-1)^{k_d-l'_d+s}, 
\end{equation}
where the generalized combination number $C_n^k$ is defined as
\begin{equation}
    \label{eq:comb_number}
    C_n^k=\left\{
    \begin{array}{ll}
        \frac{n!}{k!(n-k)!}, & 0\leqslant k\leqslant n,  \\
        0, & k>n \text{\quad or \quad } k<0.  
    \end{array}
    \right.
\end{equation}
Besides, the term $\ga_{\kappa}^{\bj}$ in \eqref{eq:theorem_A} is given by
\begin{equation}
    \label{eq:gamma}
    \begin{split}
    \ga_{\kappa}^{\bj}=\int_{\mR^3}\int_{S^2}
    &\left[H_{\bj}\left({{\frac{\bg'}{\sqrt2}}}\right)-H_{\bj}
    \left({{\frac{\bg}{\sqrt2}}}\right)\right] 
    \Hk\left({{\frac{\bg}{\sqrt2}}}\right)B(|\bg|,\sigma)
    \omega\left({{\frac{\bg}{\sqrt2}}}\right) \rd \sigma \rd \bg,
    \end{split}
\end{equation}
where from \eqref{eq:post-velo}, it holds for $\bg'$ that 
\begin{equation}
    \label{eq:coe_g}
    \bg'\triangleq \bv'-\sv'=\frac{1-e}{2}\bg+\frac{1+e}{2}|\bg|\sigma.
\end{equation}
\end{theorem}
The proof of Theorem \ref{thm:step1} is similar to \cite[Theorem 1]{Approximation2019}. For the completeness of this work, we provide it in App. \ref{app:thm1}. 
Unlike the classical case, $|\bg'|$ does not equal $|\bg|$ in the inelastic model. Therefore, $\ga_{\kappa}^{\bj}$ could not be further simplified as done in \cite{Approximation2019}. However, for special collision kernels such as the VHS kernel, \eqref{eq:gamma} could be calculated exactly, which will be discussed in the next section.

\subsection{Simplification of VHS model}
\label{sec:VHS}
For the VHS kernel \eqref{eq:VHS}, which does not depend on the collision parameter $\sigma$, the coefficient $\gamma_{\kappa}^{\bj}$ can be calculated exactly. We will begin with two lemmas. 
\begin{lemma}
\label{thm:sp_int}
Assuming $\sigma = (\sigma_1, \sigma_2, \sigma_3)$ is a unit vector, and $\kappa = (\kappa_1, \kappa_2, \kappa_3) \in \bbN^3$, then
    \begin{equation}
        \label{eq:int_sigma}
      \mS(\kappa) \triangleq  \int_{S^2}\si_1^{\kappa_1}\si_2^{\kappa_2}\si_3^{\kappa_3}\rd \sigma=
        \left\{
            \begin{array}{ll}
            4\pi\frac{(\kappa-1)!!}{(|\kappa|+1)!!}, & 2| \kappa,  \\
        0, & \text{otherwise},
            \end{array}
        \right.
    \end{equation}
where $(-1)!!$ is regarded as $1$, $\kappa !! = \kappa_1 !! \kappa_2!!\kappa_3 !!$, and  $2|\kappa$ means that all the components of $\kappa$ are even.
\end{lemma}
\begin{proof}[Proof of Lemma \ref{thm:sp_int}]
The proof can be completed with a spherical coordinate transformation.
\end{proof}

\begin{lemma}
\label{thm:int_Her}
For the Hermite polynomial $H_{\alpha}(\bv)$ and the weight function $\omega(\bv)$ defined in \eqref{eq:short}, it holds that 
\begin{equation}
    \label{eq:int_poly}
    \begin{split}
        \mathcal{V}(\kappa, \alpha, \mu) &\triangleq  \int_{\mR^3} \bv^{\kappa} \Ha(\bv)|\bv|^{\mu}\omega(\bv)\rd \bv\\
        &=(2\pi)^{-\frac32}\sum_{\sk{\bj\pq\alpha \\ 2|(\alpha-\bj)}}
        \mC(\alpha, \bj)2^{\frac{1+\mu+|\bj|+|\kappa|}2} \Gamma\left(\frac{3+\mu+|\bj|+|\kappa|}2\right)\mS(\bj+\kappa),
    \end{split}
\end{equation}
where $\alpha, \kappa \in \bbN^3$ and $\mu\in \bbR^+$. $\Gamma(\cdot)$ denotes the Gamma function and 
\begin{equation}
    \label{eq:coe_mC}
    \mC(\alpha, \beta) = \mC(\alpha_1, \beta_1)\mC(\alpha_2, \beta_2)\mC(\alpha_3, \beta_3)
\end{equation}
with $\mC(n,k), n, k\in\bbN$ being the coefficient of $x^k$ in $H_n(x)$.
\end{lemma} 
The proof of Lemma \ref{thm:int_Her} is provided in App. \ref{app:VHS}. 

Defining the coefficients $D(\alpha, \beta, \mu)$ and $\psi(\alpha, \beta, \mu)$ as 
\begin{align}
    \label{eq:coe_D}
    &D(\alpha, \beta, \mu)=\int_{\mR^3}\int_{S^2}\Ha(\bg')H_{\beta}(\bg)|\bg|^{\mu}\omega(\bg)\rd \sigma\rd \bg, \\
    \label{eq:coe_psi}
    &\psi(\alpha, \beta, \mu)=\int_{\mR^3}\int_{S^2}\Ha(\bg)H_{\beta}(\bg)|\bg|^{\mu}\omega(\bg)\rd \sigma\rd \bg,
\end{align}
then the following theorem can be established for the VHS model as 
\begin{theorem}
For the VHS kernel $B(|\bg|, \sigma) = C |\bg|^{2(1-\varpi)}$, $\gamma_{\kappa}^{\bj}$ in \eqref{eq:gamma} can be calculated exactly as 
\label{thm:VHS}
\begin{equation}
    \label{eq:cal_gamma}
    \gamma_{\kappa}^{\bj} = C 2^{\frac{5}{2}-\varpi}\left[D\Big(\bj, \kappa, 2(1-\varpi)\Big)-\psi\Big(\bj, \kappa, 2(1-\varpi)\Big)\right].     
\end{equation}
\end{theorem}
\begin{proof}[Proof of Theorem \ref{thm:VHS}]
Applying the change of variables $\bg\rightarrow \sqrt2 \bg$ and using \eqref{eq:coe_D} and \eqref{eq:coe_psi}, one can simplify \eqref{eq:gamma} to \eqref{eq:cal_gamma}.
\end{proof}
For now, the remaining task is to calculate $D(\alpha, \beta, \mu)$ and $\psi(\alpha, \beta, \mu)$ in \eqref{eq:coe_D} and \eqref{eq:coe_psi}. The result is proposed in the proposition below (detailed proof is given in App. \ref{app:VHS}). 
\begin{proposition}
\label{thm:coe_D_psi}
The coefficients $D(\alpha, \beta, \mu)$ and $\psi(\alpha, \beta, \mu)$ in \eqref{eq:coe_D} and \eqref{eq:coe_psi} can be calculated exactly as 
\begin{equation}
    \label{eq:lemma_D}
    \begin{split}
    D(\alpha, \beta, \mu)
    =\sum_{\sk{\lambda\pq\alpha \\ 2|(\alpha-\lambda)}}\mC(\al,\la)
   \sum_{\sk{\kappa\pq\la \\ 2|(\la-\kappa)}} C_{\lambda}^{\kappa}
         \left(\frac{1-e}{2}\right)^{|\kappa|}\left(\frac{1+e}{2}\right)
    ^{|\la|-|\kappa|} \mS(\lambda-\kappa) \mV(\kappa, \beta, |\lambda|-|\kappa|+\mu),
    \end{split} 
\end{equation}
and
\begin{equation}
    \label{eq:lemma_psi}
    \psi(\alpha, \beta, \mu)=4\pi \sum_{\sk{\lambda\pq\alpha \\ 2|(\alpha-\lambda)}}\mC(\al,\la)\mV(\lambda, \beta, \mu),
\end{equation}
where 
\begin{equation}
    \label{eq:com_num}
    C_{\lambda}^{\kappa} = C_{\la_1}^{\ka_1}C_{\la_2}^{\ka_2}C_{\la_3}^{\ka_3}. 
\end{equation}
The right-hand side of \eqref{eq:com_num} represents the combination number defined in \eqref{eq:comb_number}.
\end{proposition}
For Maxwell molecules ($\varpi=1$), since the collision kernel $B$ does not depend on $\bg$, we have the following proposition regarding the special sparsity of $\Aalk$. The proof is provided in App. \ref{app:VHS}. 
\begin{proposition}
    \label{thm:Maxwell}
    For the Maxwell molecules, it holds for the coefficients $A\alk$ that 
    $A\alk=0$ when $|\alpha|<|\lambda|+|\kappa|$.
\end{proposition}

\begin{table}[ht]
    \centering
    \def\arraystretch{1.5}
    {\footnotesize
    \begin{tabular}{llll}
    \hline 
    Coefficients & Formula & Used in & Computational cost \\
    \hline 
        $\mS(\kappa)$ & \eqref{eq:int_sigma} & \eqref{eq:int_poly}, \eqref{eq:lemma_D} & $\mO(M^3)$ \\
        $\mV(\kappa, \alpha, \mu)$ & \eqref{eq:int_poly} & \eqref{eq:lemma_D}, \eqref{eq:lemma_psi} & $\mO(M^{10})$ \\
       $D(\alpha, \beta, \mu)$ & \eqref{eq:coe_D} & \eqref{eq:cal_gamma} & $\mO(M^{9})$ \\
        $\psi(\alpha, \beta, \mu)$ & \eqref{eq:coe_psi} & \eqref{eq:cal_gamma} & $\mO(M^{9})$ \\
        $\ga_{\kappa}^{\bj}$ & \eqref{eq:gamma} & {\eqref{eq:theorem_A}} & $\mO(M^{6})$ \\
        $\ca{d}$ & \eqref{eq:ca} & {\eqref{eq:theorem_A}} & $\mO(M^4)$ \\
        $A\alk$ & \eqref{eq:theorem_A} & {\eqref{eq:simp1_Q}} & $\mO(M^{12})$ \\
    \hline
    \end{tabular}
    }
    \caption{The formulas and computational costs for obtaining $A\alk$ and related coefficients in the VHS model.}
    \label{table:cost1}
\end{table}


Consequently, the eight-dimensional integral in \eqref{eq:Aalk} is reduced to merely a series of summations for the VHS model. The computational cost for all related coefficients is listed in Tab. \ref{table:cost1}. It can be observed that the computational cost for calculating all $\Aalk$ is $\mO(M^{12})$, but this is not a major issue as $\Aalk$ can be pre-computed offline and stored for the simulation. 

Nevertheless, it is still computationally expensive to solve the inelastic Boltzmann equation directly using \eqref{eq:simp2_Q}. The memory required to store $\Aalk$ is $\mO(M^9)$ \cite{Approximation2019}, which is too large for practical applications. Moreover, the computational cost for each collision term is also $\mO(M^9)$, and it becomes even larger in spatially inhomogeneous tests, making it unacceptable for large values of $M$. Thus, following the strategy in \cite{Approximation2019, ZhichengHu2019, multi2022}, we adopt a special design for the numerical algorithm to reduce the computational cost. This will be introduced in the next section.

\section{Numerical scheme}
\label{sec:numerical}
In this section, we introduce the numerical scheme to solve the moment equations \eqref{eq:moment}. The Strang-splitting approach is utilized to split the moment equation into a convection step and a collision step. Specifically, the numerical scheme for the convection step is proposed in Sec. \ref{sec:con}, and the specially designed method to solve the collision step is discussed in Sec. \ref{sec:scheme_coll}.

For convenience, we first consider the numerical scheme for spatially one-dimensional spatial cases, where we have
\begin{equation}
    \label{eq:explain1D}
    \pd{\cdot}{x_2}=\pd{\cdot}{x_3}=0.
\end{equation}
Therefore, the Boltzmann equation is split into 
\begin{itemize}
    \item Convection step 
    \begin{equation}
        \label{eq:con}
        \pd{f}{t} + v_1 \pd{f}{x_1} = 0.
    \end{equation}
    \item Collision step 
    \begin{equation}
        \label{eq:col_step}
        \pd{f}{t} = \frac{1}{\Kn}\mQ[f,f](\bv).
    \end{equation}
\end{itemize}

\subsection{Convection step}
\label{sec:con}
Before introducing the numerical scheme to solve the convection step, we need to choose the expansion center $[\ou, \oT]$ in the expansion \eqref{eq:Her-expan}. Following \cite{ZhichengHu2019, multi2022},  we choose a spatially and temporally constant $[\ou, \oT]$ for the convection step, and the value depends on the specific problem. Let $\bdf^{[\ou, \oT]}$ be a column vector with all $f_{\alpha}^{[\ou, \oT]}, |\alpha| \leqslant M$ as its components. Thus, the moment equations \eqref{eq:moment} of the convection step can be rewritten in matrix-vector form as 
\begin{equation}
    \label{eq:vec_moment}
    \pd{\bdf^{[\ou, \oT]}}{t}+\bA_1^{[\ou, \oT]}\pd{\bdf^{[\ou, \oT]}}{x_1}= 0,
\end{equation}
where $\bA_1^{[\ou, \oT]}$ is a constant matrix and can be diagonalized. We refer the readers to \cite{ZhichengHu2019} for more details. 

Next, we propose the numerical scheme to solve the convection term. Suppose a spatial domain $\Omega\subset\bbR$ is discretized by a uniform grid with cell size $\Delta x$ and cell centers $\{x_j\}$. We denote $\left(\bdf_j^{[\ou, \oT]}\right)^n$ as the approximation of the average of $\bdf^{[\ou, \oT]}$ in the $j$th grid cell $[x_j-\frac12\Delta x, x_j+\frac12\Delta x]$ at time $t^n$.  The finite element method is used to solve the convection part, and the system can be solved using the forward-Euler method with a time step size $\Delta t$ as follows:
\begin{equation}
    \label{eq:scheme_con}
    \begin{split}
        &\left(\bdf^{[\ou, \oT]}\right)_j^{n+1, \ast}=\left(\bdf^{[\ou, \oT]}\right)_j^n-\frac{\Delta t}{\Delta x}\left(\bF_{j+1/2}^n-\bF_{j-1/2}^n\right), \\
        \end{split}
\end{equation}
where $\bF_{j+1/2}^n$ is the numerical flux chosen according to the HLL scheme \cite{HLL}
\begin{equation}
    \label{eq:HLL_flux}
    \bF^n_{j+1/2}=
    \left\{
    \begin{array}{ll}
    \bA_1\bdf^{n,L}_{j+1/2},& \la^L\geqslant 0, \\
    \frac{\la^R\bA_1\bdf^{n,L}_{j+1/2}-
    \la^L\bA_1\bdf^{n,R}_{j+1/2}+\la^R\la^L
    \left(\bdf^{n,R}_{j+1/2}-\bdf^{n,L}_{j+1/2}\right)}
    {\la^R-\la^L},& \la^L<0<\la^R, \\[2mm]
    \bA_1\bdf^{n,R}_{j+1/2},& \la^R\leqslant 0.
    \end{array}
    \right.
\end{equation}
Here, the superscript $[\ou, \oT]$ on $\bdf$ is omitted for simplicity in \eqref{eq:HLL_flux}. 
The characteristic velocities $\lambda^L$ and $\lambda^R$ in \eqref{eq:HLL_flux} are defined by
\begin{equation}
    \la^L=\overline{u}_1 - C_{M+1}\sqrt{\bT},\quad \la^R=\overline{u}_1 + C_{M+1}\sqrt{\bT},
\end{equation} 
where they represent the minimum and maximum eigenvalues of $\bA_1$, and $C_{M+1}$ is the largest root of the standard Hermite polynomial (defined in \eqref{eq:coef_Her}) of degree $M+1$. In \eqref{eq:HLL_flux}, $\bdf^{n,L}_{j+1/2}$ and $\bdf^{n,R}_{j+1/2}$ are computed using the WENO reconstruction method \cite{Weno}, and the details can be found in App. \ref{app:WENO}. Furthermore, the time step size must be chosen to satisfy the CFL condition
\begin{equation}
    \label{eq:CFL}
    {\rm CFL}\triangleq \Delta t \frac{|\overline{u}_1| + C_{M+1}\sqrt{\bT}}{\Delta x} < 1.
\end{equation}
Now we have completed the numerical scheme for the one-dimensional spatial case. This scheme can be naturally extended to three-dimensional spatial situations. 


\subsection{Collision step}
\label{sec:scheme_coll}
For the collision step, as stated before, the computational cost to compute the collision term is still quite expensive, on the order of $\mO(M^9)$. Therefore, we propose a new collision model to reduce the cost, following the idea in \cite{Approximation2019}. In this section, we will introduce the new collision model and then discuss the numerical scheme in the collision step.

\subsubsection{Building new collision model}
\label{sec:new_col}
To build the new collision model, both the quadratic collision term $\mQ[f, f]$ and a simplified collision operator $\mQ_S[f, f]$ are utilized. The quadratic operator is used to obtain the low-order terms in the new collision model, while the high-order terms are approximated with the simplified collision operator to save memory and computational cost. This approach has been successfully applied to the elastic Boltzmann equation \cite{Approximation2019}, where the BGK collision model is used as the simplified collision operator. 

First, the expansion center is chosen following the same method as in \cite{ZhichengHu2019}, where the local macroscopic velocity and temperature $[\bu, \theta]$ are utilized:
\begin{equation}
    \label{eq:Her-expan-local}
    f(t,\bx,\bv) \approx\sum_{|\alpha| \leqslant M}f^{[\bu,\theta]}_{\alpha}(t,\bx)
		\mH_{\alpha}^{[\bu, \theta]}(\bv).
\end{equation}
Following similar lines in \cite{Approximation2019, multi2022}, the new collision model is built by combining the quadratic collision term \eqref{eq:expan_coll} with a simplified operator $\mQ_S[f, f]$. Assume the simplified operator $\mQ_S[f, f]$ can be expanded as 
\begin{equation}
    \label{eq:coeff_QL}
    \mQ_S[f, f](\bv) \approx \sum_{|\alpha|\leqslant M} Q_{S, \alpha}^{[\bu, \theta]}(t, \bx)\mH_{\alpha}^{[\bu, \theta]}(\bv).  
\end{equation}
Then, under the expansion center $[\bu, \theta]$, the new collision model is built as 
\begin{equation}
    \label{eq:new_coll}
     \mQ_{\rm new}[f, f](\bv)  = \sum_{|\alpha| \leqslant M}  Q_{{\rm new}, \alpha}^{[\bu, \theta]}(t, \bx) \mH^{[\bu, \theta]}_{\alpha}(\bv),
     \end{equation}
with 
\begin{equation}
\label{eq:coe_new_coll}
    Q_{{\rm new}, \alpha}^{[\bu, \theta]}=\left\{ 
    \begin{array}{ll}
        \sum\limits_{|\lambda|, |\kappa| \leqslant M_0}\Aalk^{[\bu, \theta]} f^{[\bu, \theta]}_{\lambda}
		f^{[\bu, \theta]}_{\kappa},  & |\alpha| \leqslant M_0, \\[4mm]
        Q_{S,\alpha}^{[\bu, \theta]}(t, \bx), & M_0 < |\alpha|\leqslant M,
    \end{array}
    \right.  
\end{equation}
where $M_0<M$ represents the order of expansion coefficients derived from the quadratic collision term.
\begin{remark}
In \eqref{eq:coe_new_coll}, the coefficient $\Aalk^{[\bu, \theta]}$ is derived through $\Aalk$ using Proposition \ref{thm:change_center}. With this technique, the memory consumption has been reduced to $\mO(M_0^9)$. 

The parameter $M_0$ is problem-dependent and is always determined empirically. Based on our experience, $M_0 = 10$ is sufficient for most problems. 
\end{remark}

Now we will discuss how to choose the simplified collision model $\mQ_S[f, f]$. In the inelastic case, we cannot directly utilize the BGK collision model due to the energy dissipation. Instead, a simplified model with a similar form to \eqref{eq:linear} is used: 
\begin{equation}
    \label{eq:new-linear}
    \mQ_{S}[f](\bv)=\nu_1\rho\left[f_{G}(\bv)-f(\bv)\right]+\nu_2\rho\nabla_{\bv}\cdot\left[(\bv-\bu)f(\bv)\right].
\end{equation}
Unlike \eqref{eq:linear}, there is a factor of $\rho$ in each term to match the quadratic form of the distribution functions in the original collision model. $\nu_1$ and $\nu_2$ are constant parameters that will be discussed later. The expansion coefficients of $f_G$ can be computed with \cite{CaiPHD} 
\begin{equation}
    \label{eq:Her-ESBGK}
    f_{G,\alpha}^{[\bu, \theta]}=\left\{\begin{array}{ll}
        \rho, & \alpha=0, \\
        0, & |\alpha|=1, \\
        \frac{1-1/\Pr}{\alpha_i \rho}\sum_{k=1}^3\sigma_{ik}f_{G, \alpha-e_i-e_k}^{[\bu, \theta]}, & |\alpha|\geqslant 2, \; i\in\{1,2,3\} \text{ s.t. } \alpha_i>0.
    \end{array}    
    \right.
\end{equation}
Here, it should be noted that the last relationship in \eqref{eq:Her-ESBGK} holds for any $i\in\{1,2,3\}$ subject to $\alpha_i>0$. 
Using $f_0^{[\bu, \theta]} = \rho$, the energy loss term $\rho \nabla_{\bv} \cdot [(\bv - \bu) f(\bv)]$ is expanded as 
\begin{equation}
    \label{eq:new-linear-2}
     \rho \nabla_{\bv} \cdot [(\bv - \bu) f(\bv)] \approx \sum_{|\alpha|\leqslant M} f_0^{[\bu, \theta]}\left(|\alpha|f_{\alpha}^{[\bu, \theta]}+\sum_{d=1}^3f_{\alpha-2e_d}^{[\bu, \theta]}\right)	\mH_{\alpha}^{\bu, \theta}(\bv),
\end{equation}
where $f_{\alpha}^{[\bu, \theta]}$ is regarded as $0$ in \eqref{eq:Her-ESBGK} and \eqref{eq:new-linear-2} if $\alpha$ contains any negative index. 

Combining \eqref{eq:Her-ESBGK} and \eqref{eq:new-linear-2}, the expansion coefficient $\mQ_{S, \alpha}^{[\bu, \theta]}$ can be derived as 
\begin{equation}
    \label{eq:coe_QL_1}
     Q_{S, \alpha}^{[\bu, \theta]}  = \nu_1f_0^{[\bu, \theta]}\left(
    f_{G,\alpha}^{[\bu, \theta]} - f_{\alpha}^{[\bu, \theta]}\right) - \nu_2 f_0^{[\bu, \theta]}\left(|\alpha|f_{\alpha}^{[\bu, \theta]}+\sum_{d=1}^3f_{\alpha-2e_d}^{[\bu, \theta]}\right).
\end{equation}
Since $\nu_1$ indicates the damping rate of high-order terms, we follow the same approach as in \cite[Sec. 3.3.2]{Approximation2019} to determine this parameter.  The goal is to ensure that the high-order terms decay faster without introducing a gap in the damping rate between terms with $|\alpha|\leqslant M_0$ and $|\alpha|>M_0$. Thus, we consider $\Aalk$ as a matrix with respect to $\lambda$ and $\kappa$ for each fixed $\alpha$, and set $\nu_1$ to be the negative value of the minimum eigenvalue of $\Aalk$ for all $|\alpha|\leqslant M_0$ (i.e. the spectral radius of the damping rate in quadratic part). We refer the readers to \cite[Sec. 3.3.2]{Approximation2019} for more details. 

As for $\nu_2$, it is borrowed from the cooling rate in \cite[(2.16)]{Astillero2005} that 
\begin{equation}
    \label{eq:nu2}
    \nu_2=\frac{2}{3\sqrt{\pi}}(1-e^2).
\end{equation}


So far, we have derived the new collision model. This new collision model reduces the computational cost for each collision term to $\mO(M_0^9+M^3)$,  leading to significant improvements in efficiency. The numerical scheme to solve the collision step using this new collision model will be discussed in the next section.

\subsubsection{Numerical scheme to update the collision step}
\label{sec:num_col}
In this section, the numerical scheme will be presented for updating the collision step based on the new collision model. Using the vector symbol $\bdf$ as in \eqref{eq:vec_moment}, we rewrite the governing equation in the collision step \eqref{eq:col_step} as 
\begin{equation}
    \label{eq:moment_col}
    \pd{\bdf^{[\bu, \theta]}_{j}}{t} = \frac{1}{\Kn} \bQ_{\rm new}[\bdf^{[\bu, \theta]}, \bdf^{[\bu, \theta]}], 
\end{equation}
where $\bQ_{\rm new}[\bdf^{[\bu, \theta]}, \bdf^{[\bu, \theta]}]$ is a column vector with all $ Q_{{\rm new}, \alpha}^{[\bu, \theta]}, |\alpha|\leqslant M$ as its components. 

After the convection step at time $t^{n+1}$, $\left(\bdf^{[\ou, \oT]}\right)_j^{n+1, \ast}$ is obtained. Then we derive the expansion coefficients $\left(\bdf^{[\bu^{n+1, \ast}_j, \theta^{n+1, \ast}_j]}\right)_j^{n+1, \ast}$ under the expansion center $[\bu^{n+1, \ast}_j, \theta^{n+1, \ast}_j]$ using the projection algorithm in App. \ref{app:project}. The expansion center $[\bu^{n+1, \ast}_j, \theta^{n+1, \ast}_j]$ corresponds to the macroscopic velocity and temperature after the convection step at time $t^{n+1}$, which can be obtained from \eqref{eq:macro_f}. 

Next, the forward Euler scheme is adopted to update \eqref{eq:moment_col} as  
\begin{equation}
    \label{eq:scheme_col}
    \left(\bdf^{\ast}\right)_j^{n+1}=\left(\bdf^{\ast}\right)_j^{n+1, \ast}
    +\Delta t \bQ_{\rm new}\left[\left(\bdf^{\ast}\right)^{n+1, \ast}_j, \left(\bdf^{\ast}\right)^{n+1, \ast}_j\right],
\end{equation}
where $\bdf^{\ast}$ is a shorthand notation for $\bdf^{[\bu^{n+1, \ast}_j, \theta^{n+1, \ast}_j]}$. 

Finally, the projection algorithm in App. \ref{app:project} is utilized once again to obtain $\left(\bdf^{[\ou, \oT]}\right)_j^{n+1}$ based on $\left(\bdf^{\ast}\right)_j^{n+1}$, which completes the collision part and moves on to the next time step.

\subsection{Outline of the numerical algorithm}
The overall numerical scheme is summarized in Algorithm \ref{algo:inhomo}.
\begin{algorithm}[htbp]
    \caption{Numerical algorithm}
    \label{algo:inhomo}
    \begin{algorithmic}[1]
        \item Preparation: calculate and store $\Aalk$ in \eqref{eq:Aalk} with the algorithm in Sec. \ref{sec:method}.
        \item Set $n=0$, and choose an expansion center $[\ou, \oT]$ for the convection step. Calculate the initial value of $\left(\bdf^{[\ou, \oT]}\right)_{j}^0$.
        \item Determine the time step $\Delta t^n$ with the CFL condition \eqref{eq:CFL}.       
        \item Solve the convection step \eqref{eq:scheme_con} to obtain
        $\left(\bdf^{[\ou, \oT]}\right)_j^{n+1, \ast}$.
        \item Obtain the macroscopic velocity and temperature $\bu^{n+1, \ast}_{j}, \theta^{n+1, \ast}_{j}$ of $\left(\bdf^{[\ou, \oT]}\right)_j^{n+1, \ast}$ using \eqref{eq:macro_f}.
        \item Project $\left(\bdf^{[\ou, \oT]}\right)_j^{n+1, \ast}$ to the function space of expansion center $[\bu^{n, \ast}_{j}, \theta^{n, \ast}_{j}]$ with the algorithm in App. \ref{app:project}, and obtain $\left(\bdf^{[\bu_{j}^{n+1, \ast}, \theta_j^{n+1, \ast}]}\right)_j^{n+1, \ast}$.
        \item Update $\left(\bdf^{[\bu_{j}^{n+1, \ast}, \theta_j^{n+1, \ast}]}\right)_j^{n+1}$ with \eqref{eq:scheme_col}.
        \item Project  $\left(\bdf^{[\bu_{j}^{n+1, \ast}, \theta_j^{n+1, \ast}]}\right)_j^{n+1}$  to the function space of expansion center  $[\ou, \oT]$, and obtain  $\left(\bdf^{[\ou, \oT]}\right)_j^{n+1}$. 
        \item Let $n\leftarrow n+1$, and return to step 3.
    \end{algorithmic}
\end{algorithm}

\section{Numerical experiments}
\label{sec:experiment}
In this section, we present several numerical experiments to validate the Hermite spectral method for the inelastic Boltzmann equation. We start with two homogeneous cases, one with Maxwell molecules and the other with the hard sphere (HS) collision kernel. Next, two one-dimensional spatial problems and a two-dimensional spatial problem will be tested with the HS collision kernel. 
\subsection{Homogeneous experiments}
\label{eq:sec:homo_exp}
We begin by studying two homogeneous problems: the heating source problem and Haff's cooling law. In the heating source problem, we use the Maxwell model, while in Haff's cooling law, we utilize the HS model. 
\subsubsection{Heating source problem}
\label{sec:heat}    
The heating source problem was first introduced in \cite{Noije1998}, and similar studies can also be found in \cite{Hu2016, Wu2015}. The governing equation is given by 
\begin{equation}
    \label{eq:HeatBoltz}
    \pd{f}{t}-\eps\Delta_v f=\frac{1}{\Kn}\mQ[f,f](\bv),
\end{equation}
where the second term represents the effect of the heating source with the diffusion coefficient $\eps \ll 1$. In this test, the Maxwell model (i.e. $\varpi=1$ in \eqref{eq:VHS}) is utilized, and the Knudsen number is chosen such that $\frac{1}{\Kn}B = \frac{1}{4\pi}$. 

By taking $\phi = 1$ and $\phi=\bv$ in \eqref{eq:weak_1}, one can derive the conservation of mass and momentum as 
\begin{equation}
    \label{eq:heating_con}
    \rho \equiv \rho_0, \qquad \bu = \bu_0, 
\end{equation}
where $\rho_0$ and $\bu_0$ represent the initial density and macroscopic velocity, respectively.
Without loss of generality, we suppose $\rho_0 \equiv 1$ and $\bu_0 \equiv \bz$. 

By multiplying $\phi = |\bv|^2/3$ on both sides of \eqref{eq:HeatBoltz} and using the weak form \eqref{eq:weak_1}, one can derive the governing equation of the temperature $\theta$ as \cite{Hu2016} 
\begin{equation}
    \label{eq:temp_evolve}
    \begin{split}
     \pd{\theta}{t}-2\eps 
    &=-\frac{1-e^2}{4}\theta,
    \end{split}
\end{equation}
where the exact solution is given by
\begin{equation}
    \label{eq:solu_temp}
    \theta(t)=\left(\theta(0)-\frac{8\eps}{1-e^2}\right)\exp\left(-\frac{1-e^2}{4}t
    \right)+\frac{8\eps}{1-e^2}.
\end{equation}
For the numerical simulation, the expansion center is chosen as $[\ou,\bT]=[\bz,1]$. Then, the moment system of \eqref{eq:HeatBoltz} can be derived as 
\begin{equation}
    \label{eq:HeatMoment}
    \frac{\rd f_{\alpha}}{\rd t}-\eps\sum_{d=1}^3f_{\alpha-2e_d}=Q_{\rm new, \alpha}, \qquad |\alpha| \leqslant M,
\end{equation}
where $Q_{\rm new, \alpha}$ is obtained from \eqref{eq:coe_new_coll}.

\begin{figure}[!hptb]
  \centering
  \subfloat[$e=0$\label{fig:Heat1}]
  {\includegraphics[width=0.23\textwidth, height=0.18\textwidth,
    clip]{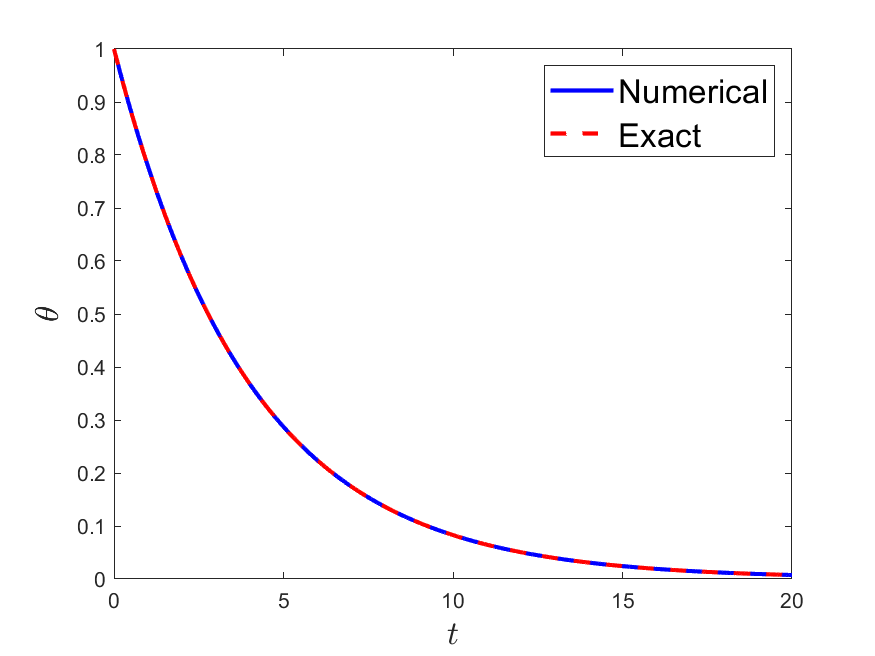}
    }  \hfill
  \subfloat[$e=0.2$]
  {\includegraphics[width=0.23\textwidth, height=0.18\textwidth,
    clip]{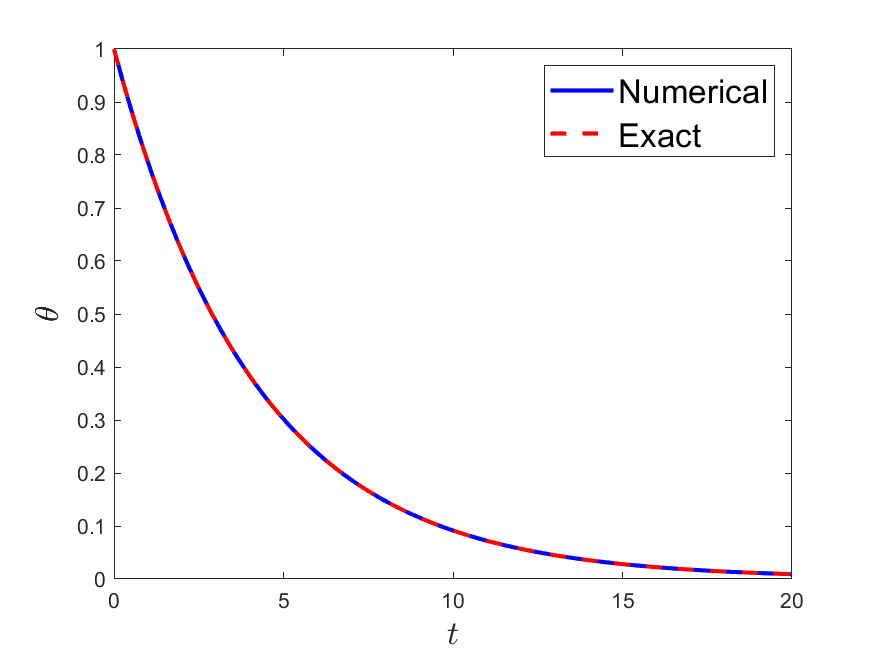} }\hfill
  \subfloat[$e=0.5$]
  {\includegraphics[width=0.23\textwidth, height=0.18\textwidth,
    clip]{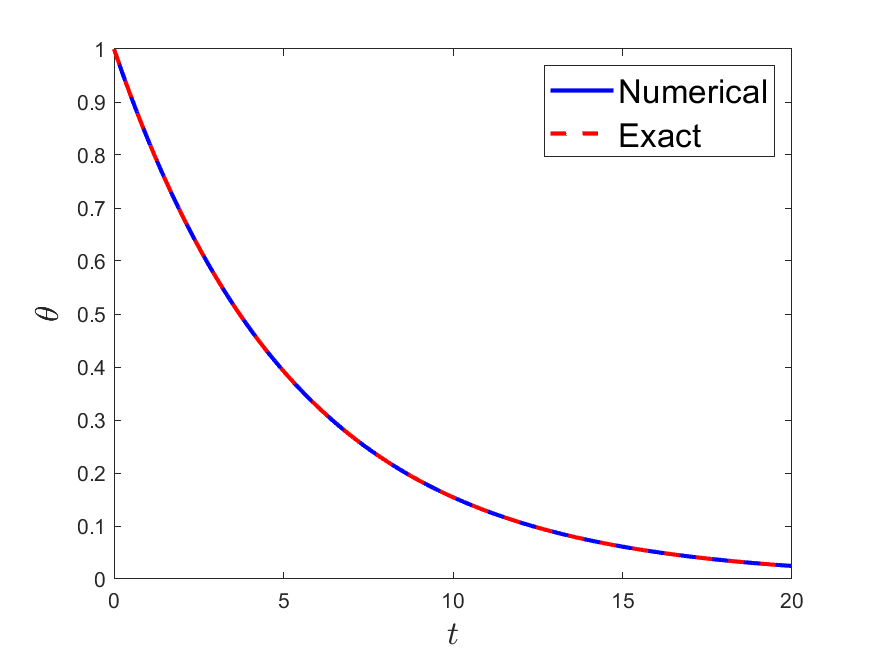}}\hfill
  \subfloat[$e=0.8$]
  {\includegraphics[width=0.23\textwidth, height=0.18\textwidth,
    clip]{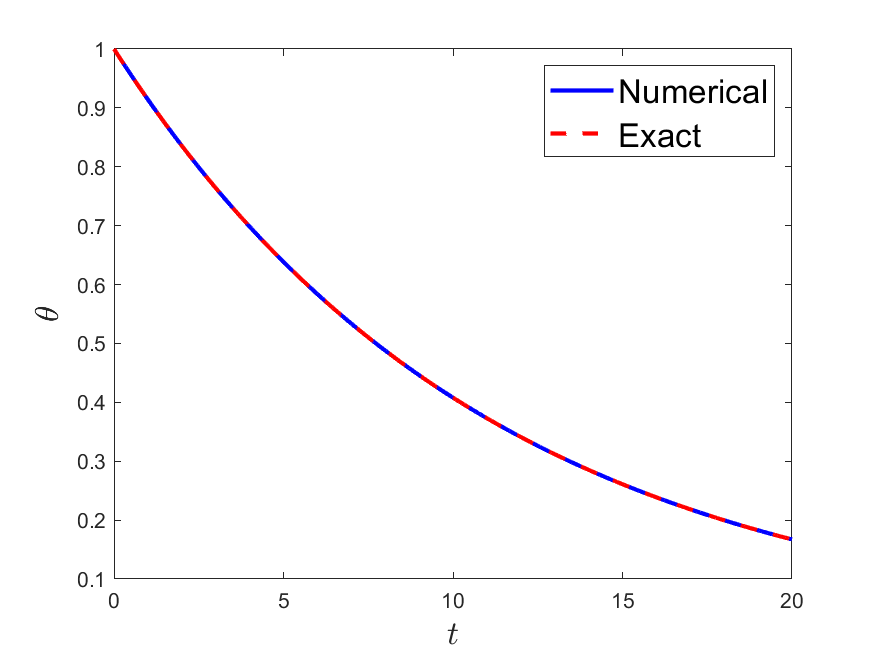}}\hfill
  \caption{(Heating source problem in Sec. \ref{sec:heat}) Evolution of temperature in the heating source problem. The restitution coefficients are $e = 0, 0.2, 0.5$ and $0.8$, respectively. The blue line represents the numerical solution, while the red line represents the exact solution in \eqref{eq:solu_temp}.}
  \label{fig:Heat}
\end{figure}

Since there is no analytical solution to this heating source problem, the solution for temperature always serves as the reference solution in this numerical test. In the simulation, we set $\epsilon = 10^{-4}$, and the initial condition is 
\begin{equation}
    \label{eq:ini_ex1}
    f(0, \bv)= \omega_{[\bz, 1]}(\bv) = \frac{1}{(2\pi)^{3/2}} \exp\left(-\frac{|\bv|^2}{2}\right).
\end{equation}
Thus, the initial condition for temperature $\theta$ is $\theta(0) = 1$. Moreover, we choose $[M_0, M] = [10, 10]$ for the length of the quadratic collision and total expansion order. The restitution coefficient $e = 0, 0.2, 0.5, 0.8$ from $t = 0$ to $t = 20$ are tested, with a time step length of $\Delta t = 0.01$. The evolution of temperature $\theta$ is displayed in Fig. \ref{fig:Heat}, which shows that the numerical solution matches well with the analytical solution for the temperature. Especially in Fig. \ref{fig:Heat1}, as time increases and the temperature decreases to zero, we can still capture the evolution of the temperature accurately. 

Here, we want to emphasize that, as shown in Proposition \ref{thm:Maxwell}, the coefficients $A\alk$ can be nonzero only when 
$|\alpha|\geqslant |\kappa|+|\lambda|$. This means the higher-order moments will not affect the lower-order ones through the collision term, and the evolution of temperature \eqref{eq:temp_evolve} is precisely described in the moment system \eqref{eq:HeatMoment}. 

We compute the error of the temperature at time $t$ as 
\begin{equation}
    \label{eq:err_theta}
    \theta_{\rm err}(t)=|\theta_{\rm num}(t)-\theta_{\rm exact}(t)|, 
\end{equation}
where $\theta_{\rm num}$ is the numerical solution and $\theta_{\rm exact}$ represents the exact solution obtained from \eqref{eq:solu_temp}. The results for different $e$ at $t=10$ are provided in Tab. \ref{table:err_Heat}. It shows that as $e$ approaches zero, the error is around $10^{-10}$, while for larger values of $e$, the error decreases even further to almost $10^{-15}$. This indicates a remarkable level of accuracy in this test.

\begin{table}[ht]
    \centering
    \def\arraystretch{1.5}
    {\footnotesize
    \begin{tabular}{l|lllll}
    \hline 
      $e$  & 0 &  0.1 & 0.2 & 0.3 & 0.4 \\
    \hline
      Error & $2.088\times10^{-10}$ & $1.437\times10^{-10}$ & $6.366\times10^{-11}$ & $1.796\times10^{-11}$ & $3.165\times10^{-12}$ \\
    \hline
    $e$ & 0.5 & 0.6 & 0.7 & 0.8 & 0.9 \\
    \hline
    Error & $3.481\times10^{-13}$ & $3.492\times10^{-14}$ & $7.383\times10^{-15}$ & $1.388\times10^{-15}$ & $1.332\times10^{-15}$ \\ 
    \hline
    \end{tabular}
    }
    \caption{(Heating source problem in Sec. \ref{sec:heat}) Error of temperature $\theta_{\rm err}$ for different $e$ at $t=10$.}
    \label{table:err_Heat}
\end{table}

\begin{remark}
    From \eqref{eq:macro_f}, it can observed that $\theta$ is expressed in terms of moments up to order $M=2$. Therefore, it is possible to obtain the evolution of temperature $\theta$ using the moment system \eqref{eq:HeatMoment} with $[M_0, M] = [2, 2]$, which further reduces the computational cost without sacrificing accuracy.
\end{remark}

\subsubsection{Haff's cooling law}
\label{sec:haff}
In this section, we numerically observe Haff's cooling law, which was first proposed by Haff in \cite{Haff1983}. The governing equation for Haff's cooling law is the same as the heating source problem as \eqref{eq:HeatBoltz} with $\epsilon = 0$. Haff's law states that for a gas composed of inelastic hard spheres, the temperature in the spatially homogeneous problem evolves as 
\begin{equation}
    \label{eq:Haff}
    \theta(t)\approx  \frac{\theta(0)}{(1+\gamma_0 t)^2}.
\end{equation}
Unlike for Maxwell molecules, the decay speed here is $\mO(t^{-2})$. Here, $\gamma_0$ is a positive constant depending on the value of $e$. We refer \cite{Hu2019, Filbet2013} for more details of this numerical test. 

\begin{figure}[!hptb]
  \centering
  \subfloat[$e=0$, $\gamma_0=0.380$ \label{fig:Haff1}]
  {\includegraphics[width=0.23\textwidth, height=0.18\textwidth,
    clip]{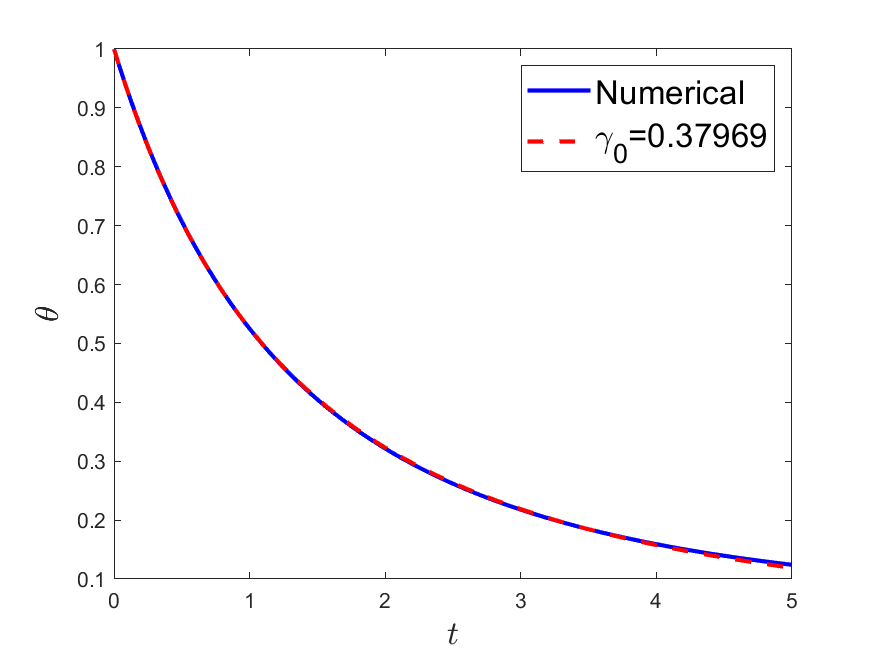}}\hfill
  \subfloat[$e=0.2$, $\gamma_0=0.364$]
  {\includegraphics[width=0.23\textwidth, height=0.18\textwidth,
    clip]{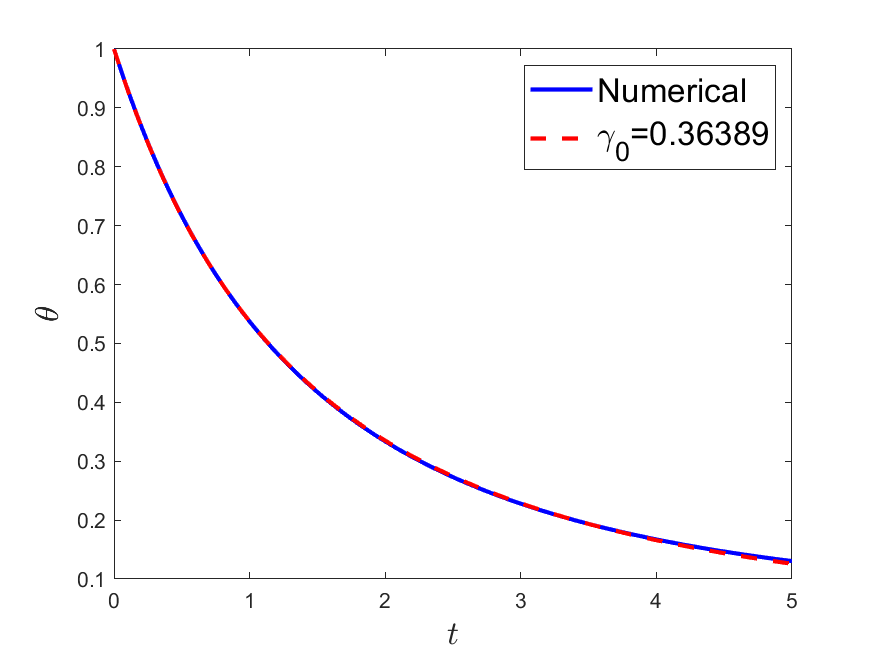}}\quad
  \subfloat[$e=0.5$, $\gamma_0=0.283$]
  {\includegraphics[width=0.23\textwidth, height=0.18\textwidth,
    clip]{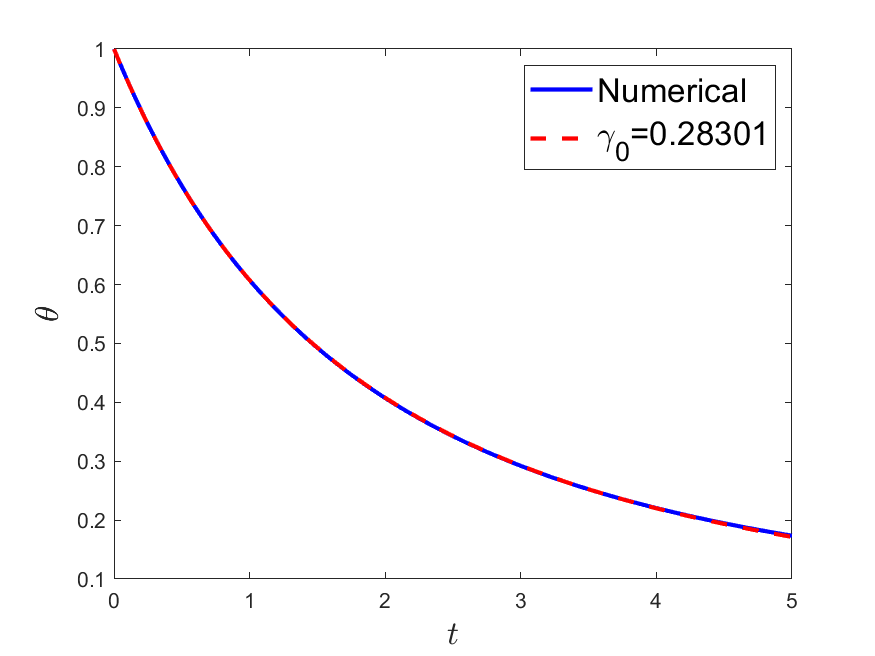}}\hfill
  \subfloat[$e=0.8$, $\gamma_0=0.135$]
  {\includegraphics[width=0.23\textwidth, height=0.18\textwidth,
    clip]{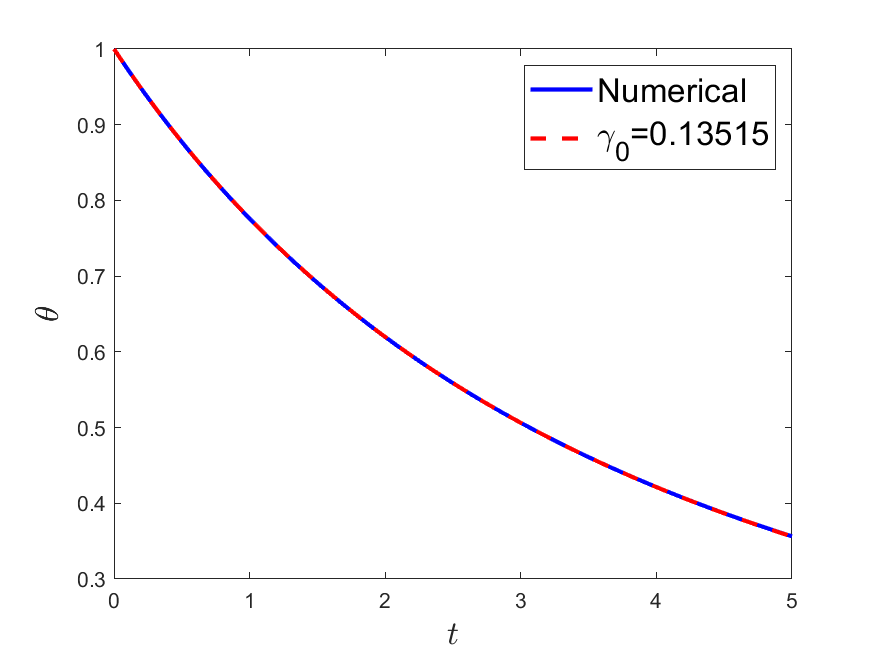}}\hfill
  \caption{(Haff's cooling law in Sec. \ref{sec:haff})
  Evolution of temperature in Haff's cooling law. The restitution coefficients are $e = 0, 0.2, 0.5$ and $0.8$, respectively. The blue line represents the numerical solution, while the red line represents the reference solution.}
  \label{fig:Haff}
\end{figure}

The nondimensionalized HS collision model has the form  
\begin{equation}
    \label{eq:HS}
    B = \frac{1}{4\sqrt{2}\pi} |\bg|.
\end{equation}
We adopt the same initial condition as in \eqref{eq:ini_ex1}, and set the Knudsen number to be $\Kn=1/\sqrt{2}$. The length of the quadratic collision and the total expansion order are chosen as $[M_0, M] = [10, 40]$, and the time step length is $\Delta t = 0.01$. The evolution of the temperature with $e = 0, 0.2, 0.5, 0.8$ from $t = 0$ to $t = 5$ is shown in Fig. \ref{fig:Haff}, where the reference solution is obtained by estimating $\gamma_0$ in \eqref{eq:Haff} using a least square fitting. From Fig. \ref{fig:Haff}, it can be clearly observed that even in the case $e = 0$, the numerical solution matches well with the reference solution. 


\subsection{Inhomogeneous experiments}
\label{sec:inhomo}
In this section, we study two one-dimensional spatial problems: Couette flow and Fourier heat transfer, as well as a two-dimensional spatial periodic diffusion problem. For all these tests, the hard sphere (HS) model is utilized as the collision model. 
The Knudsen number is calculated with
\begin{equation}
    \label{eq:express_Kn}
    \Kn=\frac{m_0}{\sqrt{2}\pi \rho_0 d_{\rm ref}^2 x_0},
\end{equation}
where the parameters in \eqref{eq:express_Kn} correspond to the nondimensionalization parameters of the working gas and HS collision kernel, which are listed in Tab.  \ref{table:character}. The method of nondimensionalization is described in App. \ref{app:nondim}. The reference solutions for these tests are obtained using the DSMC method provided in \cite{Astillero2005} for the HS collision kernel.

\subsubsection{Couette flow}
\label{sec:couette}
\begin{figure}[!hptb]
    \centering
    \subfloat[Density, $\rho$ (kg$\cdot$ m$^{-3}$)]
    {\includegraphics[width=0.45\textwidth, height=0.36\textwidth,
      clip]{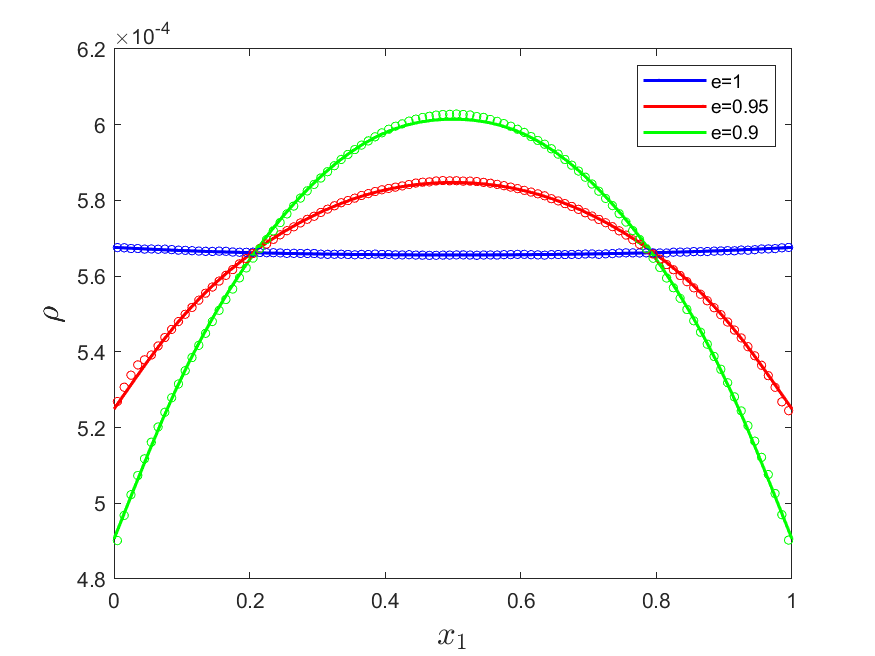}}\hfill
      \subfloat[$y$-component velocity, $u_2$ (m/s)]
    {\includegraphics[width=0.45\textwidth, height=0.36\textwidth,
      clip]{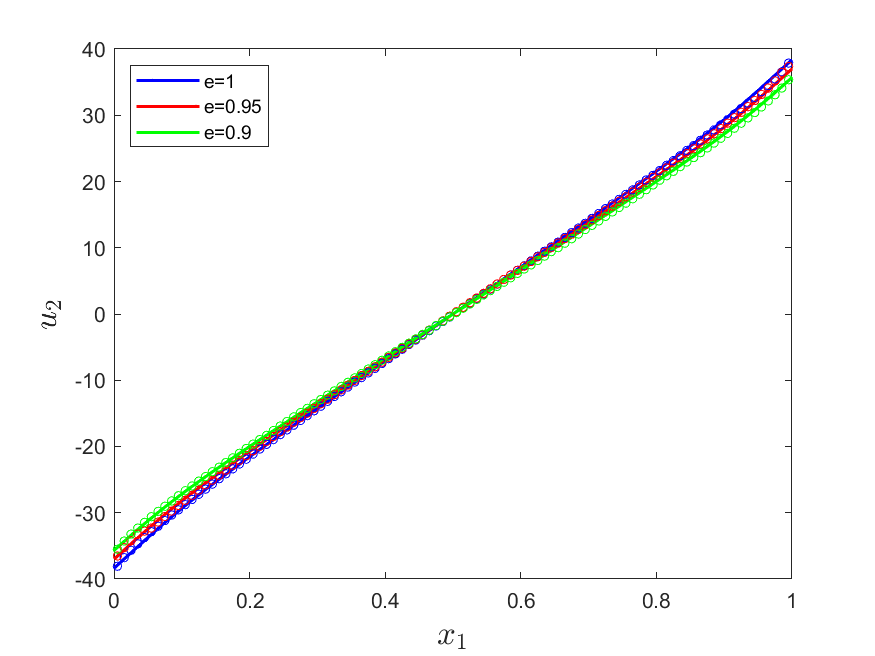}}\\
      \subfloat[Temperature, $\theta$ (K)]
    {\includegraphics[width=0.45\textwidth, height=0.36\textwidth,
      clip]{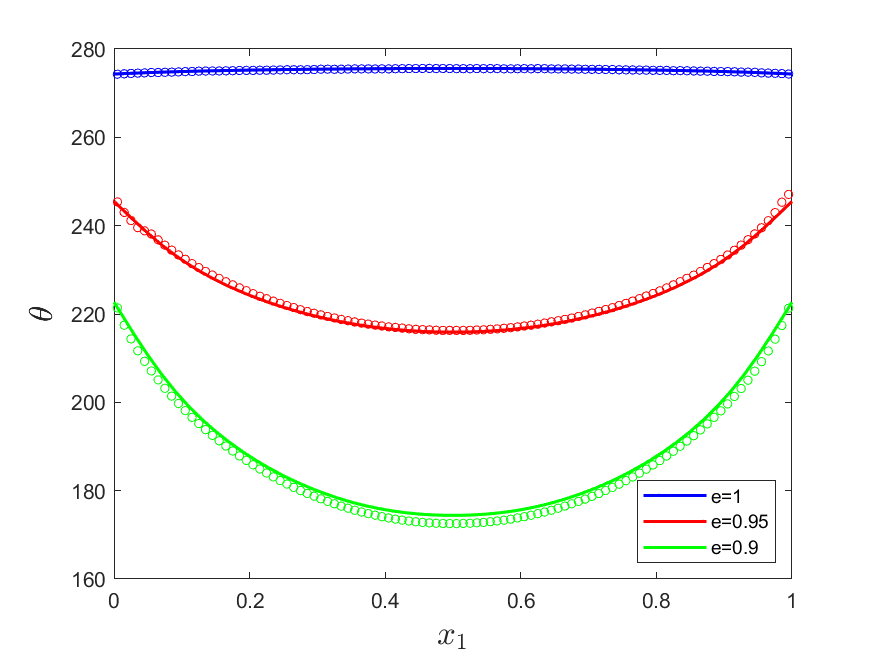}}\hfill
      \subfloat[Heat flux, $q_1$ (kg $\cdot$ s$^{-3}$) ]
    {\includegraphics[width=0.45\textwidth, height=0.36\textwidth,
      clip]{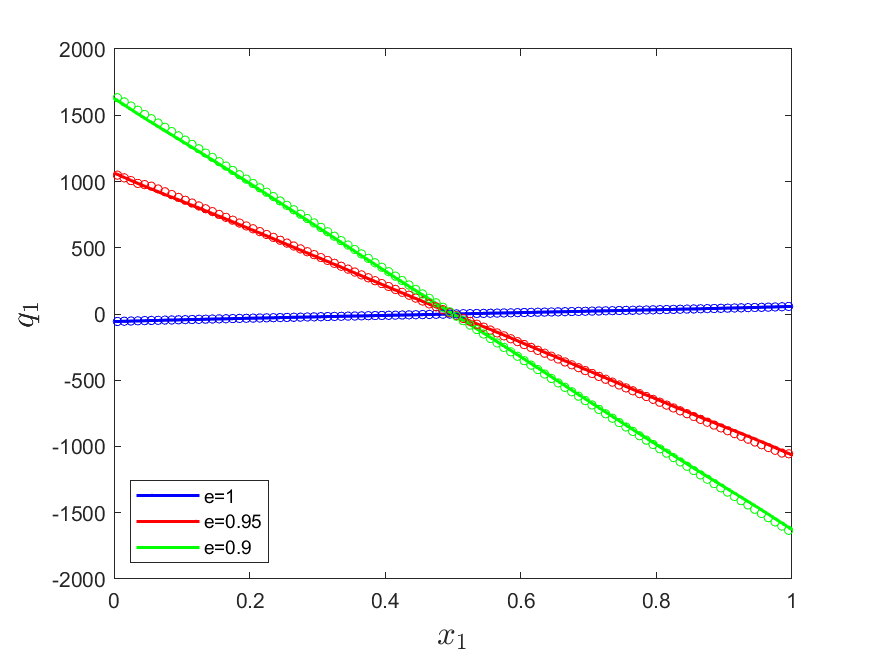}}\hfill
    \caption{(Couette flow in Sec. \ref{sec:couette}): Numerical solutions of the Couette flow for $\Kn=0.2$
    with $e = 1, 0.95$ and $0.9$.  Lines correspond to numerical solutions, and symbols denote the reference solutions from DSMC.}
    \label{fig:Couette1}
\end{figure}

In this section, we consider the 1D Couette flow, which is a benchmark problem also tested in \cite{Wu2015, ZhichengHu2019}. The setup consists of two infinite parallel plates with a distance of $10^{-3}$m. Both plates are purely diffusive and have a temperature of $273$K. They move in opposite directions along the $y$-axis with speeds $\bu^w= (0, \mp50, 0)$ m/s. The initial state is set as velocity $\bu = \bz$m/s and $\theta = 273$K. Two different densities are considered: $\rho = 5.662\times 10^{-4}$kg$\cdot$m$^{-3}$ and $1.132 \times 10^{-4}$kg$\cdot$m$^{-3}$, which correspond to $\Kn = 0.2$ and $1$, respectively.

\begin{figure}[!hptb]
    \centering
      \subfloat[Density, $\rho$ (kg$\cdot$ m$^{-3}$)]
    {\includegraphics[width=0.45\textwidth, height=0.36\textwidth,
      clip]{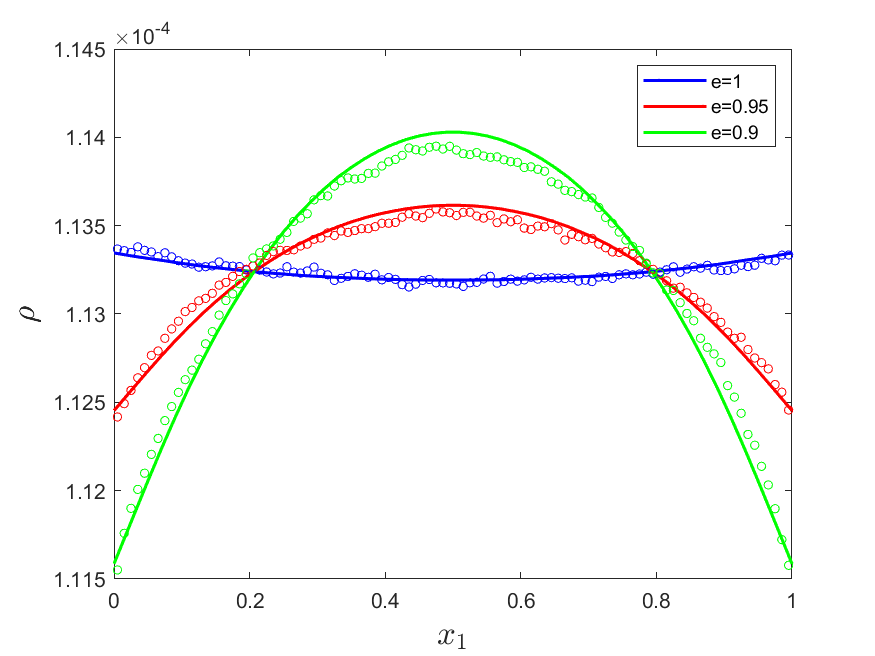}}\hfill
      \subfloat[$y$-component velocity, $u_2$ (m/s)]
    {\includegraphics[width=0.45\textwidth, height=0.36\textwidth,
      clip]{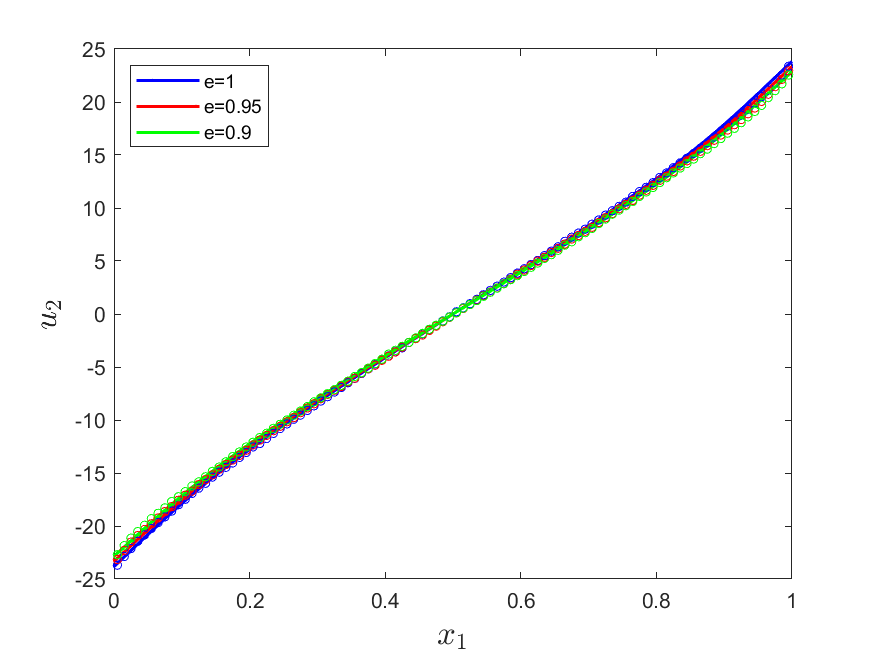}}\\
      \subfloat[Temperature, $\theta$ (K)]
    {\includegraphics[width=0.45\textwidth, height=0.36\textwidth,
      clip]{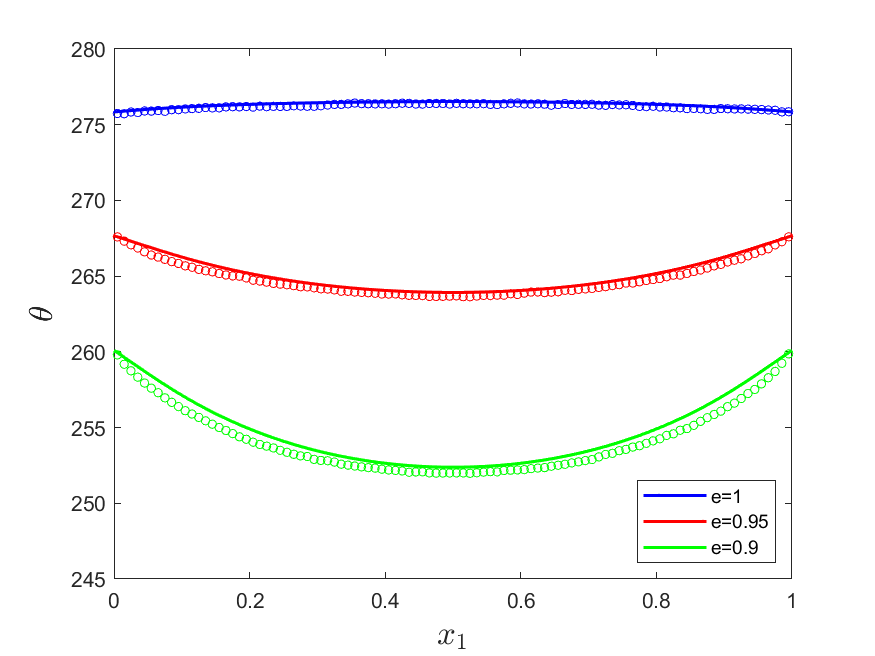}}\hfill 
      \subfloat[Heat flux, $q_1$ (kg $\cdot$ s$^{-3}$) ]
    {\includegraphics[width=0.45\textwidth, height=0.36\textwidth,
      clip]{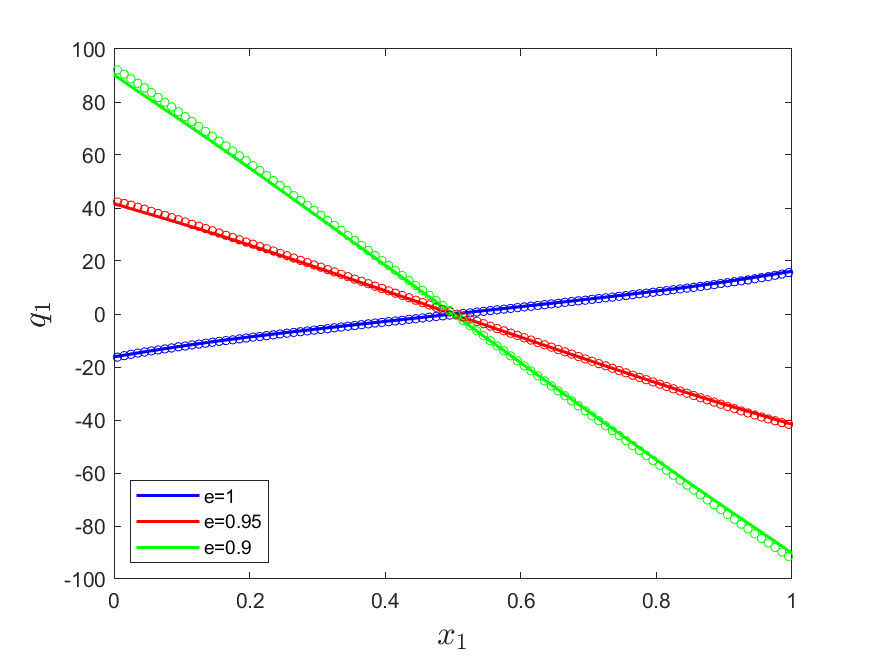}}\hfill
    \caption{(Couette flow in Sec. \ref{sec:couette}) Numerical solutions of the Couette flow for $\Kn=1.0$
    with $e = 1, 0.95$ and $0.9$.  Lines correspond to numerical solutions, and symbols denote the reference solutions from DSMC.}
    \label{fig:Couette2}
\end{figure}

In the simulation, a uniform grid with $50$ cells and WENO reconstruction are utilized for the spatial discretization, and the CFL number is set as ${\rm CFL} = 0.3$. The length for the quadratic collision term and the total expansion number are chosen as $[M_0, M] =[10, 40]$. The expansion center $[\ou, \oT]$ in the convection step is set as $[\bz, 1]$, and the restitution coefficients $e=1, 0.95$, and $0.9$ are implemented. The density $\rho$, the macroscopic velocity $u_2$ in the $y$-direction, the temperature $\theta$, and the heat flux $q_1$ in the $x$-direction at the steady state are studied. 

Numerical results for $\Kn = 0.2$ and $1$ are illustrated in Fig. \ref{fig:Couette1} and \ref{fig:Couette2}, respectively. For $\Kn = 0.2$, all the numerical solutions coincide well with the reference solutions. For the case of $\Kn = 1$, the velocity $u_2$, temperature $\theta$ and heat flux $q_1$ agree well with the reference solutions, while there is a small discrepancy in the density $\rho$, with the largest relative error being less than $1\%$. It is worth noting that there are some oscillations in the reference results, while the numerical solutions remain smooth.

\subsubsection{Fourier heat transfer}
\label{sec:fourier}
\begin{figure}[!hptb]
    \centering
    \subfloat[Density, $\rho$ (kg$\cdot$ m$^{-3}$)]
    {\includegraphics[width=0.45\textwidth, height=0.36\textwidth,
      clip]{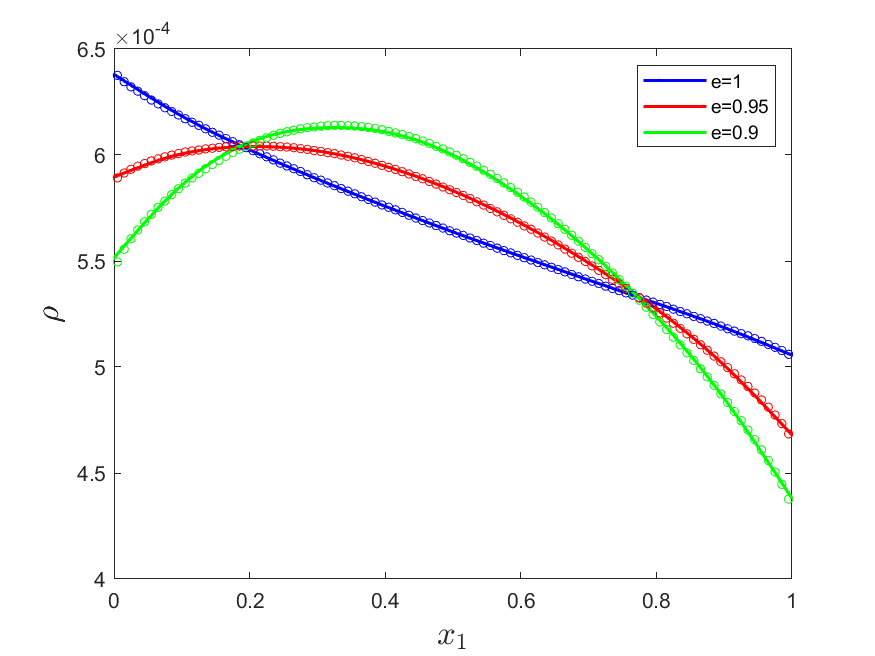}}\hfill
      \subfloat[Temperature, $\theta$ (K)]
    {\includegraphics[width=0.45\textwidth, height=0.36\textwidth,
      clip]{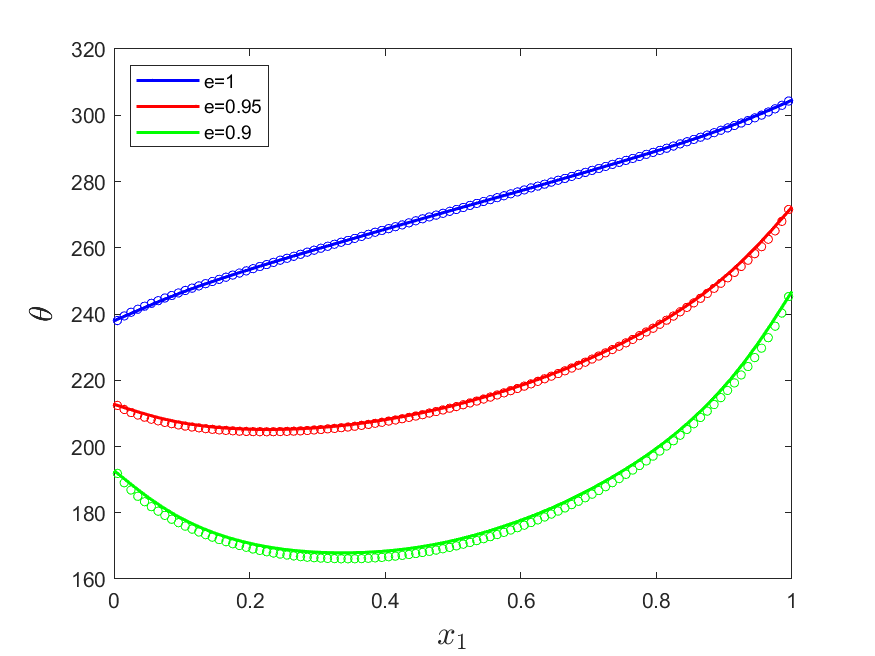}}\\
      \subfloat[Stress tensor, $\sigma_{11}$ (kg$\cdot$ m$^{-1}$ $\cdot$ s$^{-2}$)]
    {\includegraphics[width=0.45\textwidth, height=0.36\textwidth,
      clip]{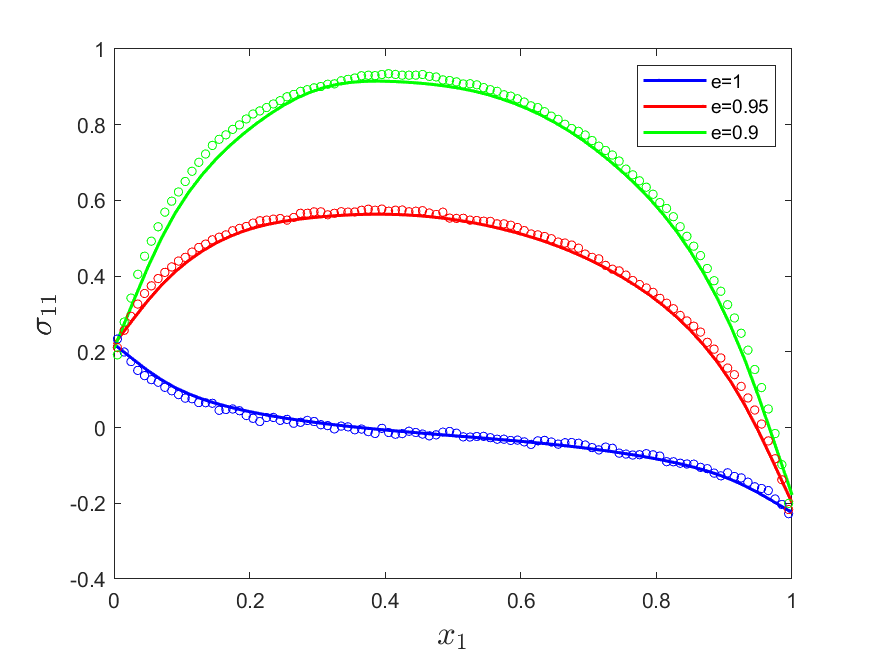}}\hfill
      \subfloat[Heat flux, $q_1$ (kg $\cdot$ s$^{-3}$) ]
    {\includegraphics[width=0.45\textwidth, height=0.36\textwidth,
      clip]{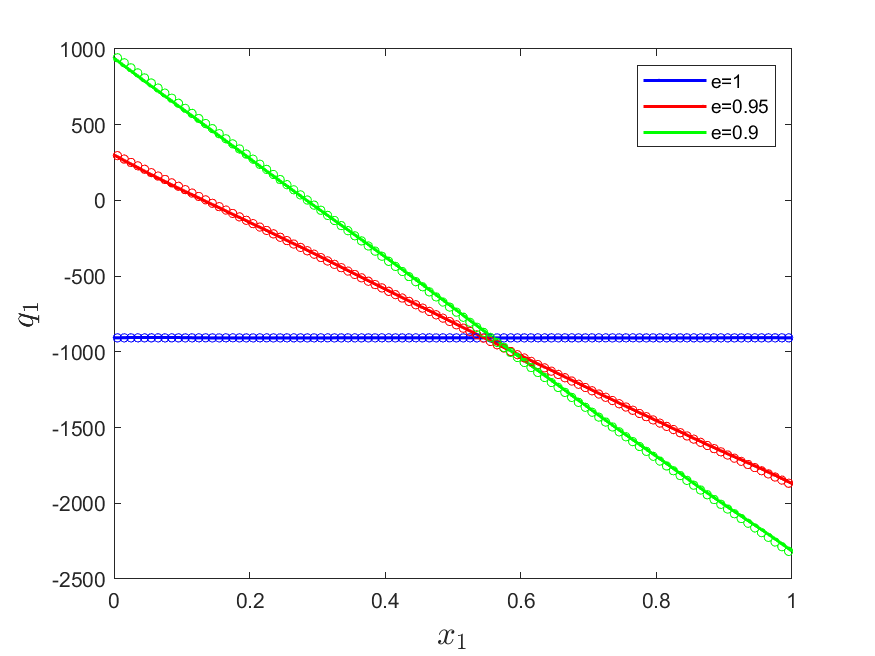}}\hfill
    \caption{(Fourier heat transfer in Sec. \ref{sec:fourier}) Numerical solutions of the Fourier heat transfer for $\Kn=0.2$ with $e = 1, 0.95$ and $0.9$. Lines correspond to numerical solutions, and symbols denote the reference solutions from DSMC.}
    \label{fig:Fourier1}
\end{figure}

Fourier heat transfer is another widely studied problem, which is also considered in \cite{Wu2015}. Similar to the Couette flow, we consider the particles between two infinitely large parallel plates. The distance between the plates is still $10^{-3}$m, and both boundaries are purely diffusive. However, in the Fourier heat transfer problem, the two plates are stationary but have different temperatures. In this case, the temperatures of the two walls are set as $\theta_l=223$K and $\theta_r=323$K. The initial conditions are set as $\bu = \bz$m/s for velocity, and $\theta = 273$K for temperature. The same two densities and values of $\Kn$ as in the Couette flow case (Sec. \ref{sec:couette}) are considered.

\begin{figure}[!hptb]
    \centering
      \subfloat[Density, $\rho$ (kg$\cdot$ m$^{-3}$)]
    {\includegraphics[width=0.45\textwidth, height=0.36\textwidth,
      clip]{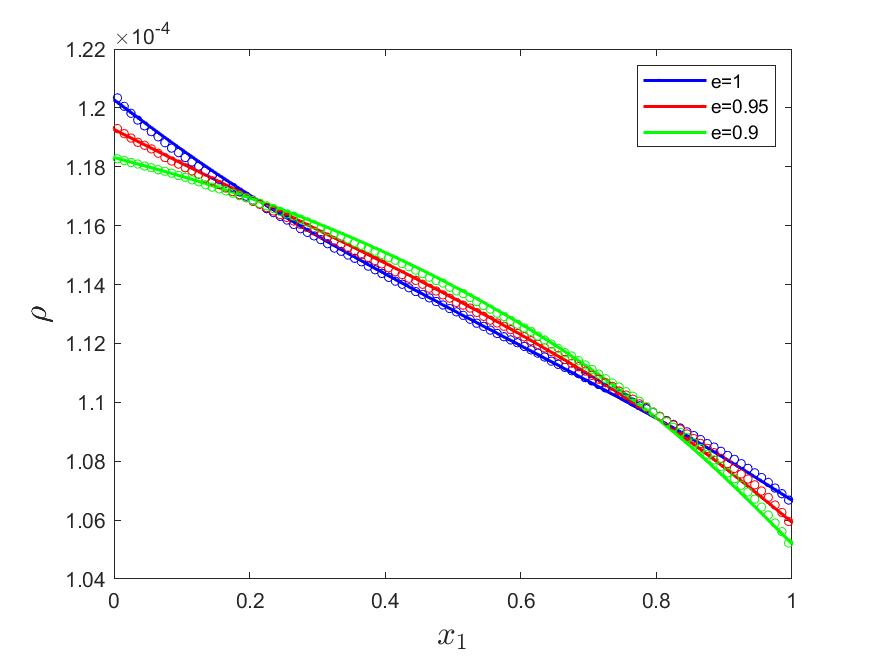}}\hfill
      \subfloat[Temperature, $\theta$ (K)]
    {\includegraphics[width=0.45\textwidth, height=0.36\textwidth,
      clip]{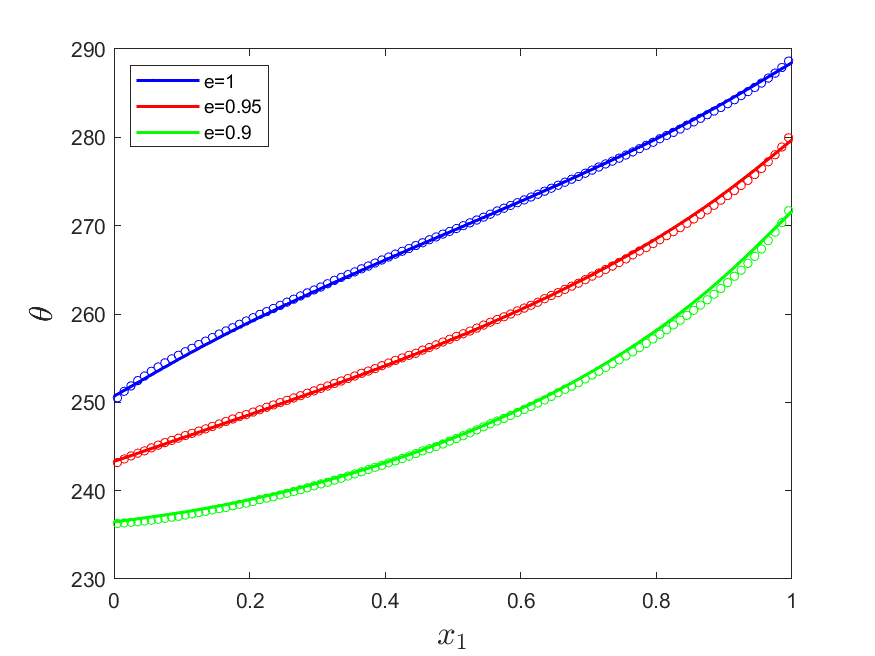}}\\
      \subfloat[Stress tensor, $\sigma_{11}$ (kg$\cdot$ m$^{-1}$ $\cdot$ s$^{-2}$)]
    {\includegraphics[width=0.45\textwidth, height=0.36\textwidth,
      clip]{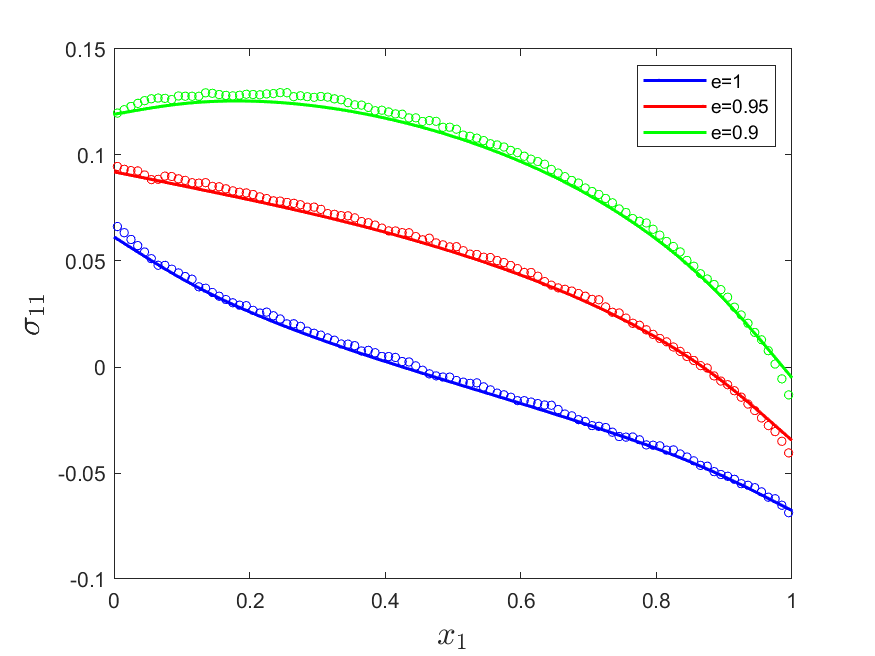}}\hfill 
      \subfloat[Heat flux, $q_1$ (kg $\cdot$ s$^{-3}$) ]
    {\includegraphics[width=0.45\textwidth, height=0.36\textwidth,
      clip]{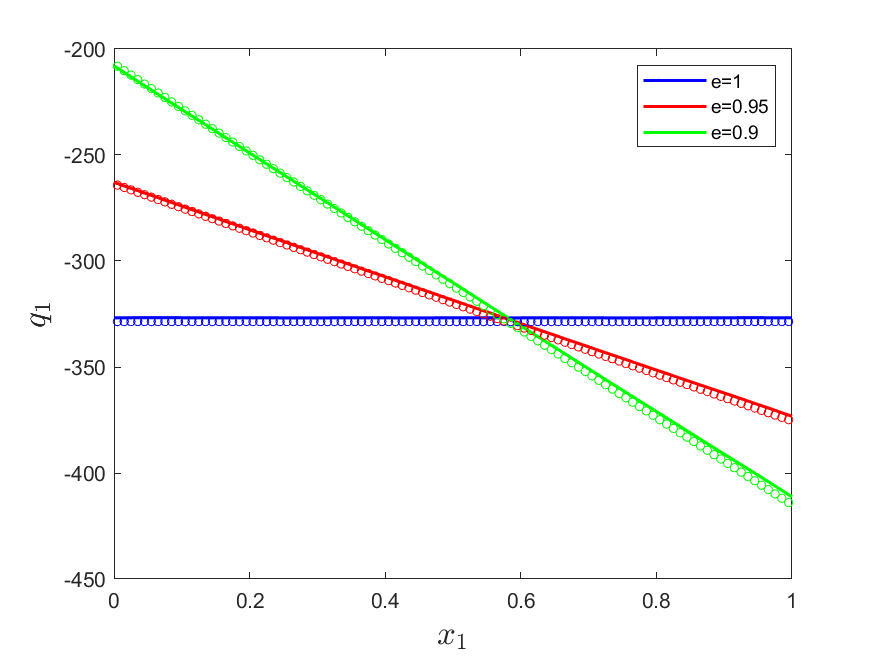}}\hfill
    \caption{(Fourier flow in Sec. \ref{sec:fourier}) Numerical solutions of the Fourier heat transfer for $\Kn=1.0$ with $e = 1, 0.95$ and $0.9$. Lines correspond to numerical solutions, and symbols denote the reference solutions from DSMC.}
    \label{fig:Fourier2}
\end{figure}

Besides, the same numerical settings such as grid, CFL number, expansion center, etc., used in the Couette flow simulation (Sec. \ref{sec:couette}), are applied here. The numerical results for $\Kn = 0.2$ and $1$ are plotted in Fig. \ref{fig:Fourier1}, and \ref{fig:Fourier2}, respectively, where the density $\rho$, temperature $\theta$, the shear stress $\sigma_{11}$, and heat flux $q_1$ at the steady state are illustrated. For both $\Kn$, the numerical solutions match well with the reference solutions, with the largest relative deviation being less than $0.5\%$ in all cases. Additionally, unlike the reference solutions by DSMC, the numerical results keep smooth.


\subsubsection{2D case: periodic diffusion}
\label{sec:diffuse}
\begin{figure}[!hptb]
  \centering
  \subfloat[$\rho$ (kg$\cdot$ m$^{-3}$), $t=0.05$]
  {\includegraphics[width=0.3\textwidth, height=0.24\textwidth,
    clip]{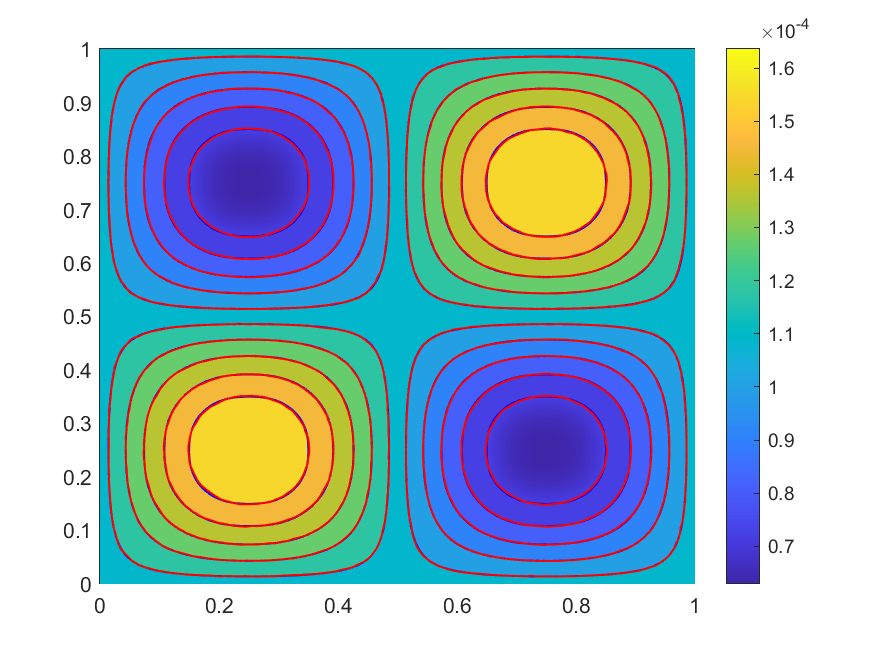}}\hfill
  \subfloat[$\theta$ (K), $t=0.05$]
  {\includegraphics[width=0.3\textwidth, height=0.24\textwidth,
    clip]{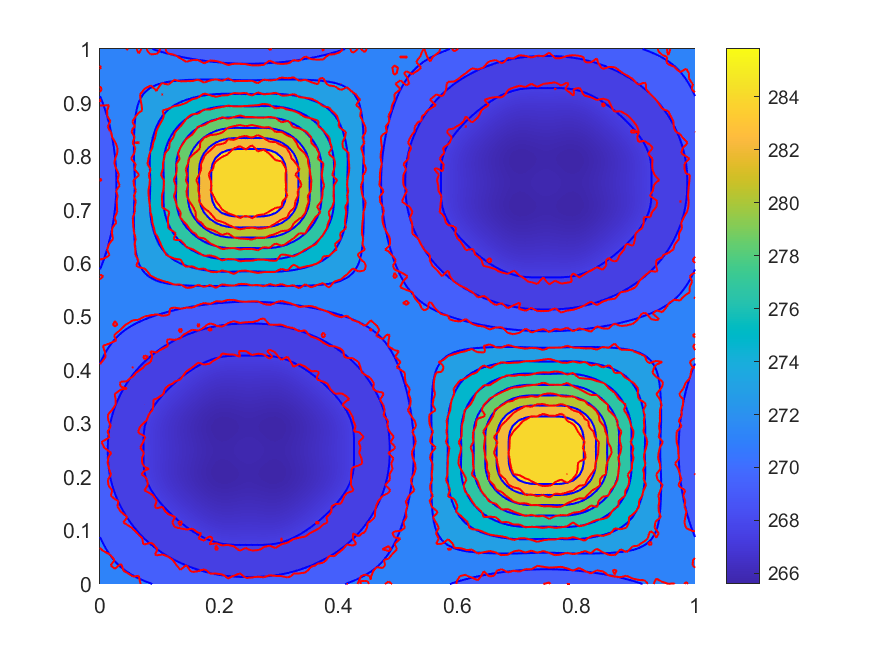}}\hfill
  \subfloat[$\sigma_{12}$ (kg$\cdot$ m$^{-1}$ $\cdot$ s$^{-2}$), $t=0.05$]
  {\includegraphics[width=0.3\textwidth, height=0.24\textwidth,
    clip]{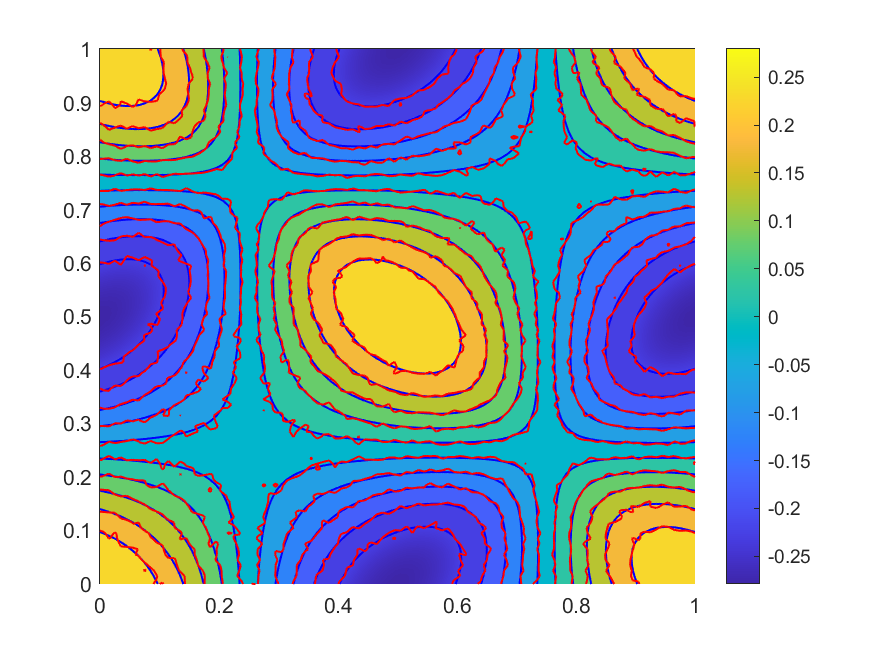}}\\
  \subfloat[$\rho$ (kg$\cdot$ m$^{-3}$), $t=0.1$]
  {\includegraphics[width=0.3\textwidth, height=0.24\textwidth,
    clip]{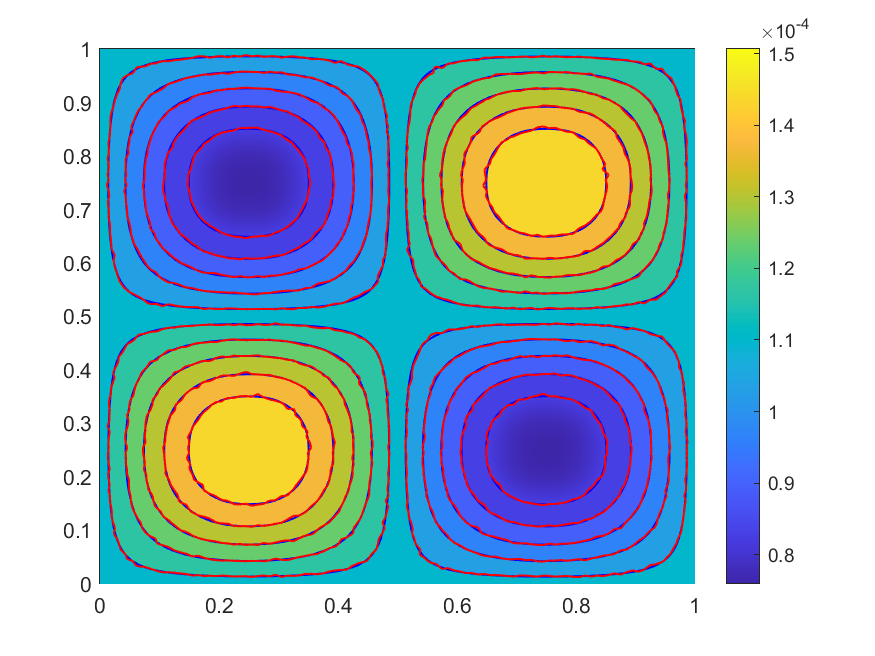}}\hfill
  \subfloat[$\theta$ (K), $t=0.1$]
  {\includegraphics[width=0.3\textwidth, height=0.24\textwidth,
    clip]{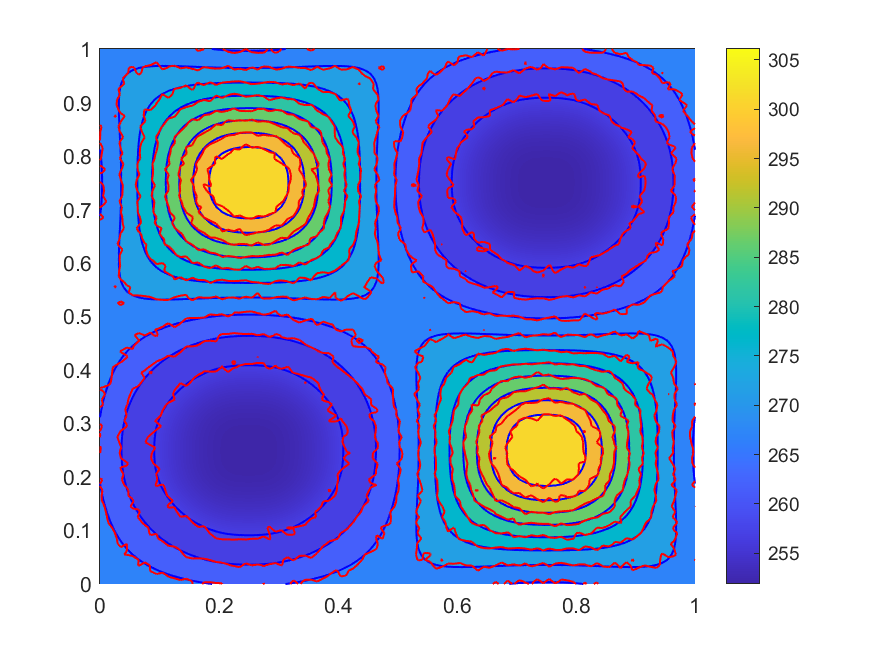}}\hfill
  \subfloat[$\sigma_{12}$ (kg$\cdot$ m$^{-1}$ $\cdot$ s$^{-2}$), $t=0.1$]
  {\includegraphics[width=0.3\textwidth, height=0.24\textwidth,
    clip]{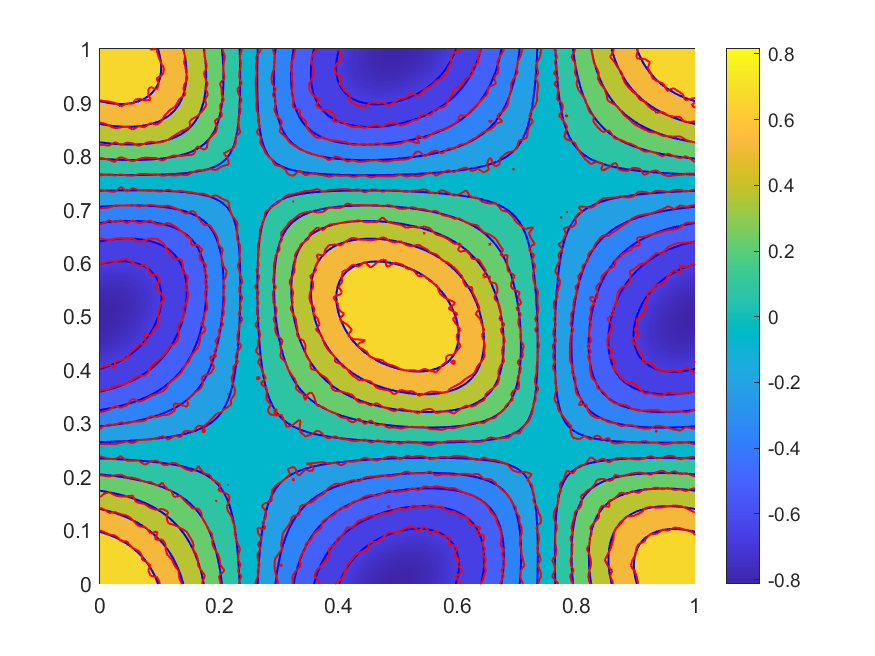}}\hfill
  \caption{(2D case: periodic diffusion in Sec. \ref{sec:diffuse}, Example 1) Solutions of the periodic diffusion for $e=0.9$ with the initial condition \eqref{eq:ini2D}. Blue contours: Numerical solutions. Red contours: Reference solutions by DSMC.}
  \label{fig:Diffuse1}
\end{figure}

In this section, we consider a two-dimensional test in the square region $\Omega=[0,L]\times [0,L]$ with periodic boundary conditions, where $L=10^{-3}$m. To validate the efficiency of the Hermite spectral method, two examples with different initial conditions are tested. 

\paragraph{Example 1}
For the first example, the initial velocity and temperature are set as $\bu=\bz$ and $\theta=273K$, respectively, throughout the entire domain. The initial density $\rho(x,y)$ is given by
{\small 
\begin{equation}
    \label{eq:ini2D}
    \rho(x,y)=1.132\times 10^{-4}\times \left[1+0.5\sin\left(2\pi\frac{x}{L}\right)\sin\left(2\pi\frac{y}{L}\right)\right] \text{kg}\cdot\text{m}^{-3}.
\end{equation}
}

In the classical case, the distribution function will diffuse to reach a global equilibrium. From a macroscopic perspective, the macroscopic variables will eventually become spatially uniform. However, in the inelastic case, due to the dissipation of total energy, the evolution of non-equilibrium macroscopic variables becomes much more complicated. Specifically, the temperature approaches zero as time increases, leading to significant challenges in the simulations.

A uniform grid with $100\times 100$ cells and the WENO reconstruction are employed for spatial discretization. With the given initial condition, the corresponding Knudsen number is $\Kn = 1$. The quadratic length and total expansion number are set as $[M_0, M] = [10, 30]$, and the expansion center for the convection step is $[\ou, \bT] = [\bz, 1]$.

 The restitution coefficient $e=0.9$ is implemented in the simulation. The macroscopic variables, including the density $\rho$, temperature $\theta$, and stress tensor $\sigma_{12}$, are studied. The numerical results are presented in Fig. \ref{fig:Diffuse1}, with reference solutions obtained from DSMC used for comparison. It can be observed that the numerical solutions agree well with the reference solutions. Furthermore, while the reference solutions exhibit some oscillations, the numerical solutions remain smooth. 

\paragraph{Example 2}
To further validate the efficiency of this method, a more complicated initial condition is considered. In this case, the variation period for the density becomes smaller, and a disturbance is introduced in the temperature as 
{\small 
\begin{equation}
    \label{eq:ini2D2}
    \begin{split}
       &\rho(x,y)=1.132\times 10^{-4}\times \left[1+0.5\sin\left(2\pi\frac{x}{L}\right)\sin\left(4\pi\frac{y}{L}\right)\right] \text{kg}\cdot\text{m}^{-3}, \\ 
       &\theta(x,y)=273\times \left[1+0.05\sin\left(2\pi\frac{x}{L}\right)\sin\left(2\pi\frac{y}{L}\right)\right] \text{K}.
    \end{split} 
\end{equation}
}

\begin{figure}[!hptb]
  \centering
  \subfloat[$\rho$ (kg$\cdot$ m$^{-3}$), $t=0.05$]
  {\includegraphics[width=0.3\textwidth, height=0.24\textwidth,
    clip]{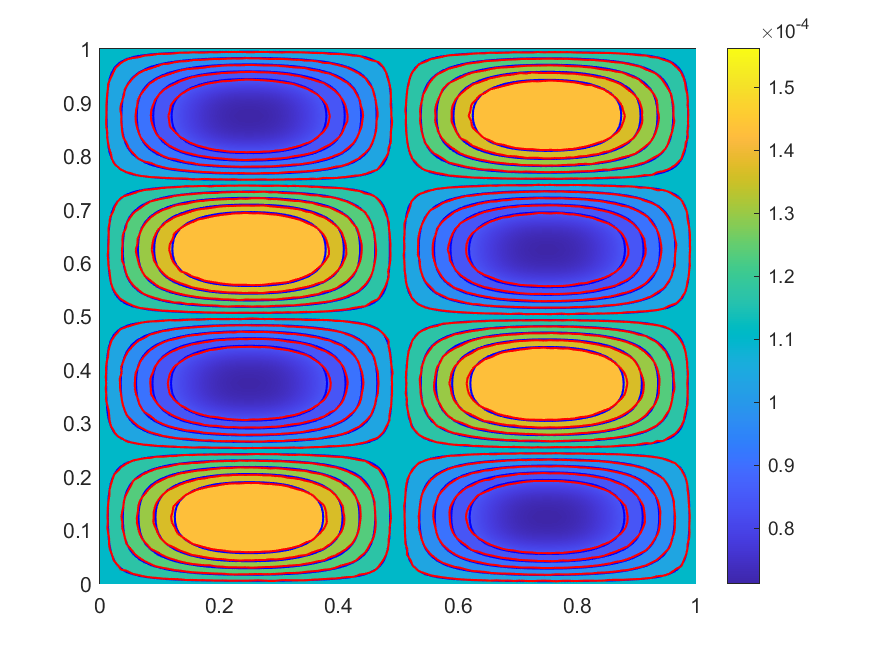}}\hfill
  \subfloat[$\theta$ (K), $t=0.05$]
  {\includegraphics[width=0.3\textwidth, height=0.24\textwidth,
    clip]{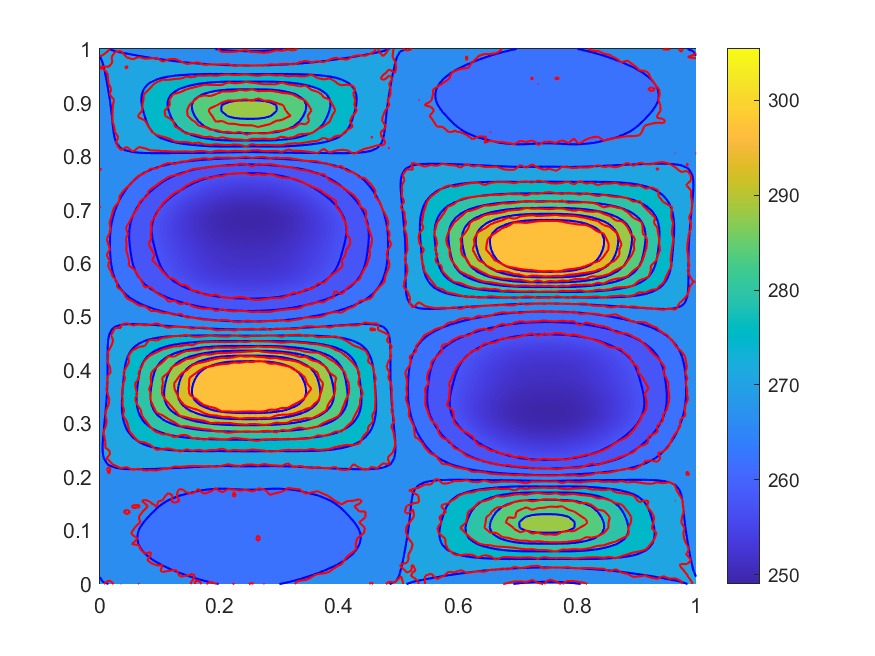}}\hfill
  \subfloat[$\sigma_{12}$ (kg$\cdot$ m$^{-1}$ $\cdot$ s$^{-2}$), $t=0.05$]
  {\includegraphics[width=0.3\textwidth, height=0.24\textwidth,
    clip]{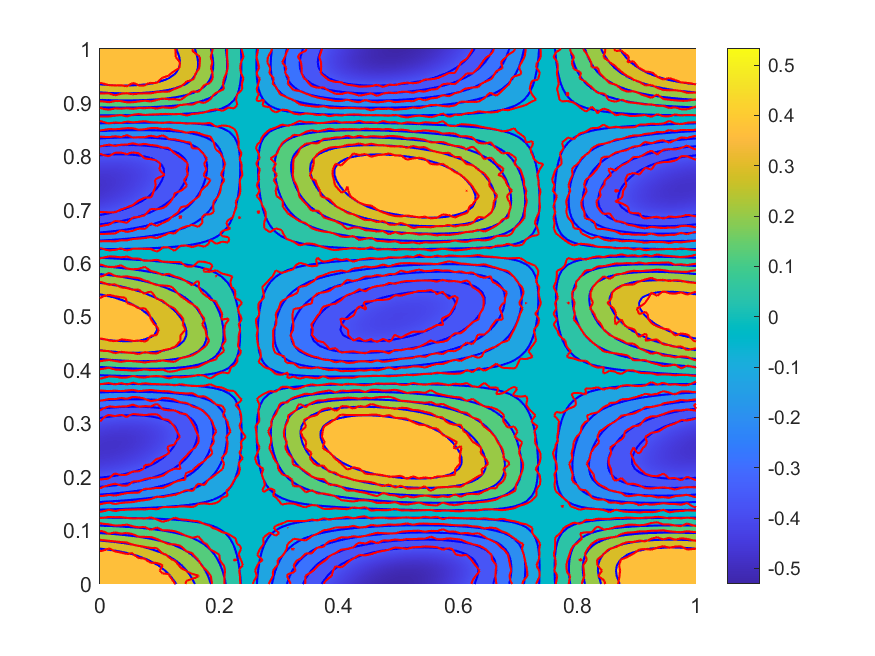}}\\
  \subfloat[$\rho$ (kg$\cdot$ m$^{-3}$), $t=0.1$]
  {\includegraphics[width=0.3\textwidth, height=0.24\textwidth,
    clip]{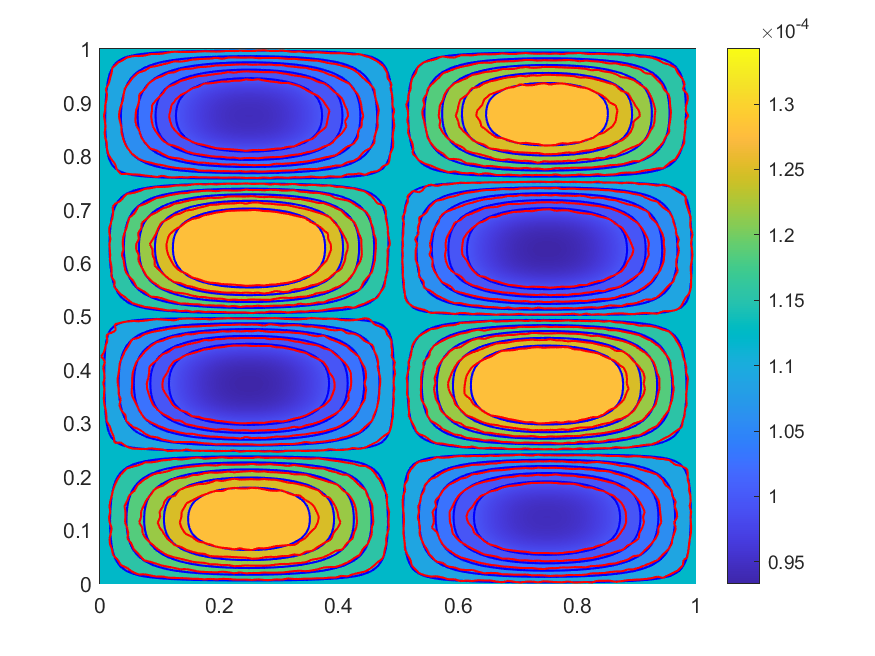}}\hfill
  \subfloat[$\theta$ (K), $t=0.1$]
  {\includegraphics[width=0.3\textwidth, height=0.24\textwidth,
    clip]{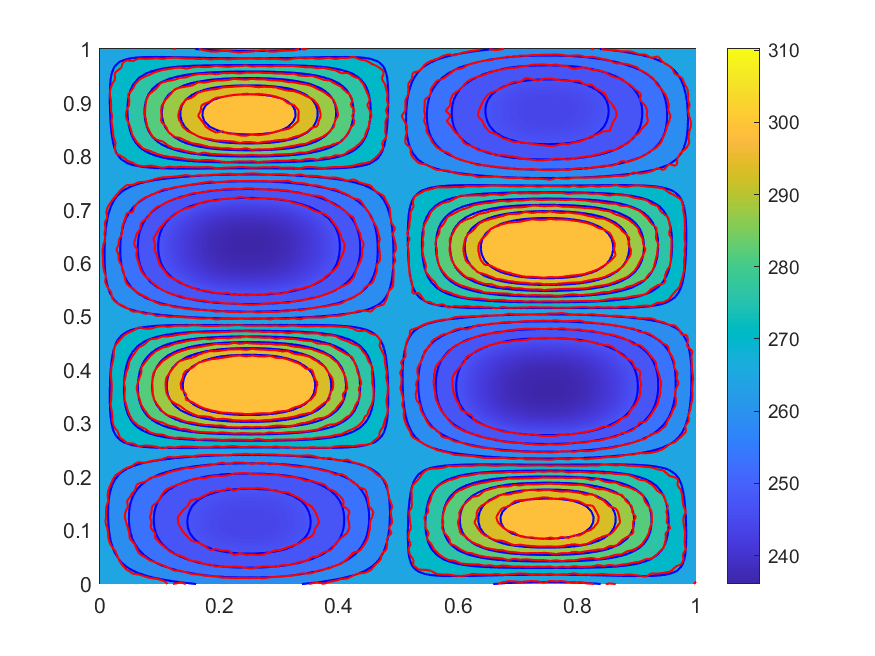}}\hfill
  \subfloat[$\sigma_{12}$ (kg$\cdot$ m$^{-1}$ $\cdot$ s$^{-2}$), $t=0.1$]
  {\includegraphics[width=0.3\textwidth, height=0.24\textwidth,
    clip]{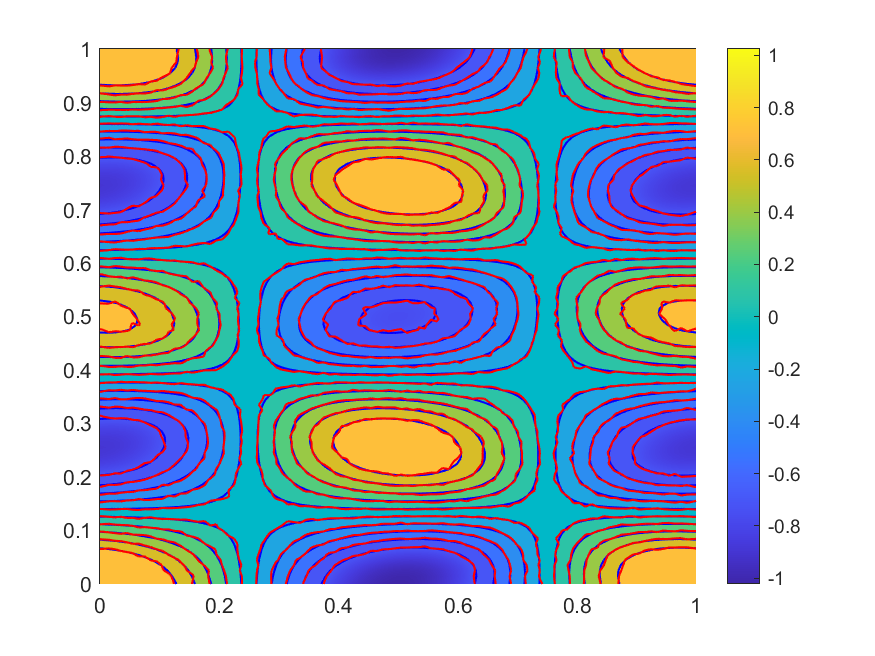}}
  \caption{(2D case: periodic diffusion in Sec. \ref{sec:diffuse}, Example 2) Solutions of the periodic diffusion for $e=0.8$ with the initial condition \eqref{eq:ini2D2}. Blue contours: Numerical solutions. Red contours: Reference solutions by DSMC.}
  \label{fig:Diffuse2}
\end{figure}

The numerical settings, such as the mesh, expansion order, etc., remain the same as Example 1, but a smaller restitution coefficient $e = 0.8$ is examined. The numerical solutions for density $\rho$, temperature $\theta$, and stress tensor $\sigma_{12}$ at $t = 0.05$ and $t = 0.1$ are shown in Fig. \ref{fig:Diffuse2}. It can be observed that even for this complex initial condition, the numerical solutions still agree well with the reference solutions. The trends of these macroscopic variables are similar to {Example 1}, while the behavior of temperature appears to be more intricate.

\begin{figure}[!hptb]
  \centering
    \subfloat[$\rho$ (kg$\cdot$ m$^{-3}$), $t=0.2$]
  {\includegraphics[width=0.3\textwidth, height=0.24\textwidth,
    clip]{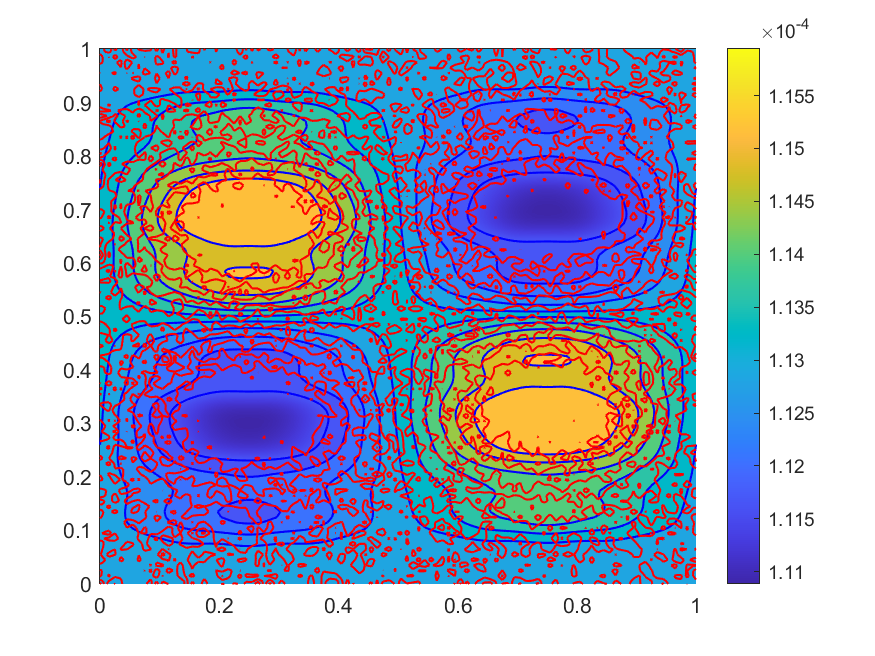}}\hfill
  \subfloat[$\theta$ (K), $t=0.2$]
  {\includegraphics[width=0.3\textwidth, height=0.24\textwidth,
    clip]{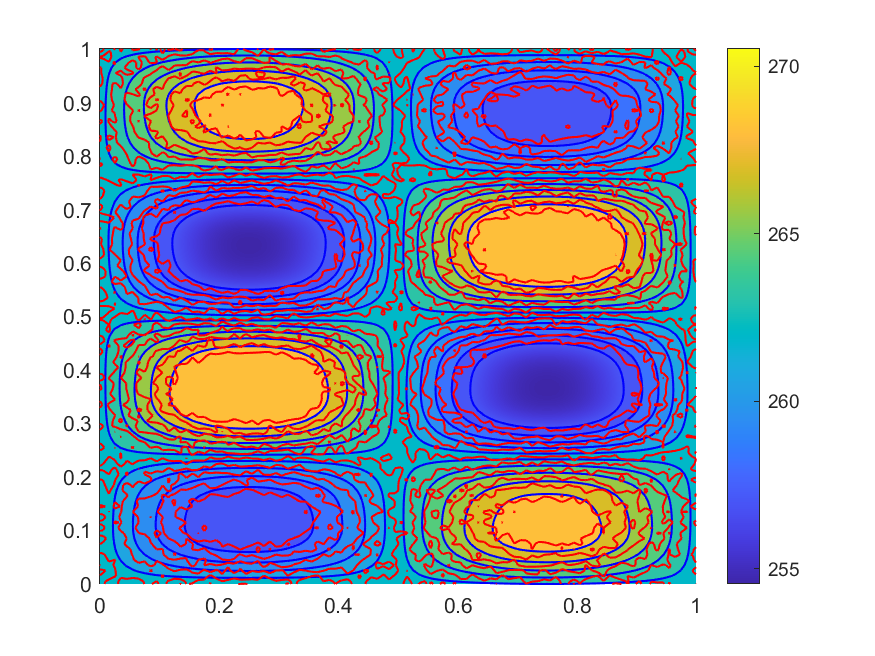}}\hfill
  \subfloat[$\sigma_{12}$ (kg$\cdot$ m$^{-1}$ $\cdot$ s$^{-2}$), $t=0.2$]
  {\includegraphics[width=0.3\textwidth, height=0.24\textwidth,
    clip]{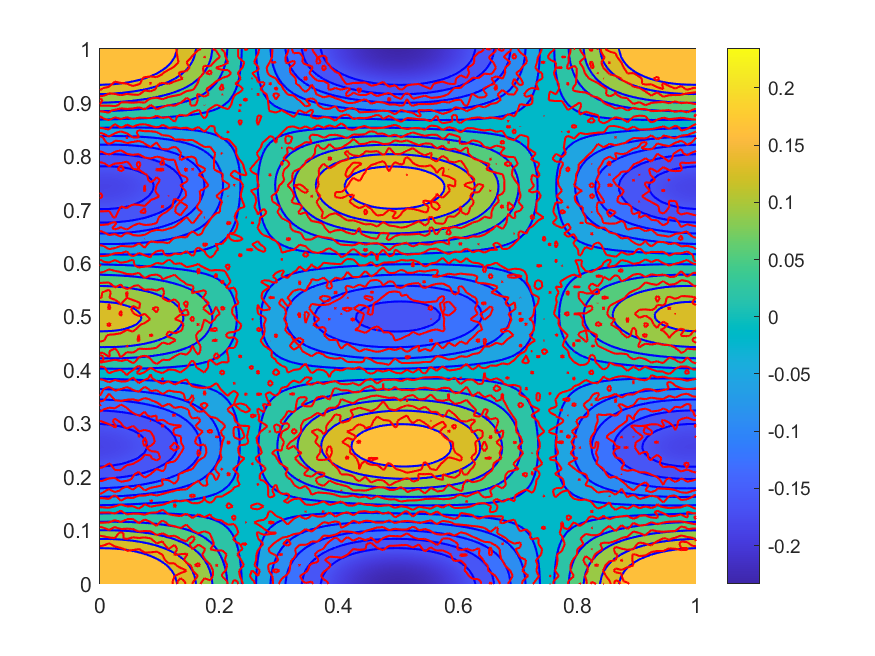}}
  \caption{(2D case: periodic diffusion in Sec. \ref{sec:diffuse}, Example 2) Solutions of the periodic diffusion for $e=0.8$ with the initial condition \eqref{eq:ini2D2} at $t = 0.2$. Blue contours: Numerical solutions. Red contours: Reference solutions by DSMC.}
  \label{fig:Diffuse3}
\end{figure}

To investigate the long-term behavior of this example, the numerical solutions at $t = 0.2$ are displayed in Fig. \ref{fig:Diffuse3}. It can be seen that the three macroscopic variables are becoming spatially uniform, while the temperature is globally decreasing. It is worth noting that the reference solutions by DSMC are filled with oscillations, which cannot capture this long-term behavior, whereas the numerical solutions of the Hermite spectral method are still smooth.

\paragraph{Efficiency test}
To quantify the efficiency of this method, we examine the computational time for both examples in the case of $t=0.1$. The simulations are performed on the CPU model Intel Xeon E5-2697A V4 @ 2.6GHz with $8$ threads utilized. The total CPU time and wall time, as well as the CPU time of each time step and grid, are provided in Tab. \ref{table:diffuse_time}.

It shows that the total time for $e = 0.9$ and $e = 0.8$ is almost the same, indicating that the restitution coefficient has a negligible effect on the computational time. Moreover, the total CPU time is almost $8$ times of the elapsed time, which indicates the excellent parallel efficiency of this Hermite spectral method. Additionally, the total degrees of freedom (DOF) in the microscopic velocity space can be calculated with 
\begin{equation}
    \label{eq:dof}
    {\rm DOF} = \frac{(M+2)(M+1)M}{6}.
\end{equation}
Hence, the total DOF in this 2D problem is $4960$ and the CPU time per DOF per grid shown in Tab. \ref{table:diffuse_time} is on the order of $\mO(10^{-6})$. These results demonstrate the high efficiency of this Hermite spectral method, making it suitable for parallel computing in large-scale problems. 

\begin{table}[!hptb]
\centering
\def\arraystretch{1.5}
{\footnotesize
\begin{tabular}{lll}
\hline
& Example 1       & Example 2          \\
\hline
$[M_0,M]$ & $[10,30]$ & $[10,30]$ \\
Restitution coefficient $e$ & 0.9 & 0.8 \\                                 
End time $t$ & 0.1 & 0.1 \\ 
\hline
Run-time data:  &  &   \\
Total CPU time $T_{\rm CPU}$ (s) & 79980 & 81454 \\ 
Elapsed time (Wall time) $T_{\rm Wall}$ (s):  & 10575.8 & 10856.3\\
Parallel efficiency & $94.53\%$ &  $93.79\%$ \\
CPU time per time step (s) & 242.36 & 246.83 \\
Degree of freedom &   $4960$ & $4960$ \\
CPU time per DOF per grid (s) & $4.9\times 10^{-6}$ & $5.0\times10^{-6}$ \\
\hline
\end{tabular}
}
\caption{(2D case: periodic diffusion in Sec. \ref{sec:diffuse}) Run-time data for the two-dimensional periodic diffusion of $t=0.1$.}
\label{table:diffuse_time}
\end{table}

\section{Conclusion}
\label{sec:conclusion}
In this paper, we have developed a numerical scheme for solving the inelastic Boltzmann equation based on the Hermite spectral method. This method demonstrates its capability to compute two-dimensional periodic model problems and accurately describe the evolution of macroscopic quantities. The expansion coefficients of the quadratic collision model are computed using the properties of the Hermite basis functions, which can be calculated exactly for the VHS model. To balance accuracy and computational cost, we introduce a new collision model that combines the quadratic collision term with a linearized collision operator.

The numerical method is validated through several benchmark problems in granular flow. Even for two-dimensional cases, the method shows excellent performance in capturing the behavior of inelastic gas flow with high efficiency.



\section*{Acknowledgements}
We thank Prof. Jingwei Hu from UW for her valuable suggestions. We thank Prof. Lei Wu from SUSTC for his help with the DSMC code. 
The work of Yanli Wang is partially supported by the National Natural Science Foundation of China (Grant No. 12171026, U2230402 and 12031013).


\begin{appendix}
\label{app:app}

\section{Proof of Theorem \ref{thm:step1}}
\label{app:thm1}
To prove Thm. \ref{thm:step1}, we first introduce the lemma below:
\begin{lemma} 
\label{lemma:tran-Hermite}
Define $\bv=\bh+\frac12\bg$, $\bw=\bh-\frac12\bg$, then 
\begin{equation}
    \label{eq:tran-Hermite}
    \Hl(\bv)\Hk(\bw)=\sum_{\ka'+\la'=\ka+\la}\ca1\ca2\ca3
    H_{\la'}(\sqrt{2}\bh)H_{\ka'}\left(\frac1{\sqrt2}\bg\right),
\end{equation}
where the coefficients $\ca{d}$ are defined in \eqref{eq:ca}. 
\end{lemma} 
The proof of Lemma \ref{lemma:tran-Hermite} can be referred to \cite[Lemma 3]{Approximation2019}. Besides, we can derive the corollary
\begin{corollary} 
\label{corollary:tran-Hermite1}
Define $\bv=\bh+\frac12\bg$, then it holds that 
\begin{equation}
    \label{tran-Hermite1}
    \Hl(\bv)=2^{-\frac{|\la|}{2}}\sum_{\ka'+\la'=\la}
    \frac{\la!}{\ka'!\la'!}
    H_{\ka'}(\sqrt{2}\bh)H_{\la'}\left(\frac1{\sqrt2}\bg\right).
\end{equation}
\end{corollary}
The proof is straightforward by letting $\ka=\bz$ in \eqref{eq:tran-Hermite}. 

Now, we can present the proof of Theorem \ref{thm:step1}.
\begin{proof}[Proof of Theorem \ref{thm:step1}]
First, we rewrite \eqref{eq:Aalk} as
\begin{equation}
    \label{eq:Aalk2}
	\begin{split}
	A\alk=&\frac{1}{\alpha!}\int_{\mR^3}\int_{\mR^3}\int_{S^2}B
	(|\bg|,\sigma)
	\Hl(\bv)\Hk(\sv)[\Ha(\bv')-\Ha(\bv)]\omega(\bv)\omega(\sv)\rd\sigma \rd\bv \rd\sv.
	\end{split}
\end{equation}
Define $\bh=\frac{\bv+\sv}2=\frac{\bv'+\sv'}2$ and note $\bg=\bv-\sv$ is the relative velocity. Besides, with $\bg'$ defined in \eqref{eq:coe_g}, we have 
\begin{equation}
    \label{eq: dd_con}
    \begin{split}
        &\bv=\bg+\frac12\bh, \qquad \sv=\bg-\frac12\bh, \qquad \bv'=\bg'+\frac12\bh,\\
        &\rd \bg \rd \bh=\rd \bv \rd \bv_{\ast}, \qquad \omega(\bv)\omega(\sv)=\omega\left(\frac1{\sqrt2}\bg\right)\omega(\sqrt2 \bh).
    \end{split}
\end{equation}
By applying Lemma \ref{lemma:tran-Hermite}, Corollary \ref{corollary:tran-Hermite1} and \eqref{eq: dd_con}, we can transform \eqref{eq:Aalk2} into an integral with respect to $\bg$ and $\bh$:
\begin{equation}
    \label{eq:Aalk3}
    \Aalk=2^{\frac{|\al|}{2}}
    \sum_{\la'+\ka'=\la+\ka}\sum_{\bi+\bj=\al}
    \frac{1}{\bi!\bj!}\ca1\ca2\ca3
    \gamma_{\ka'}^{\bj}\eta_{\la'}^{\bi},
\end{equation}
where $\gamma_{\ka'}^{\bj}$ is the integral involving $\bg$ defined in \eqref{eq:gamma}, and $\eta_{\la'}^{\bi}$ is the integral involving $\bh$ defined as
\begin{equation}
    \label{eq:eta}
    \eta_{\la'}^{\bi}=\int_{\mR^3}H_{\la'}(\sqrt{2}\bh)H_{\bi}
    (\sqrt{2}\bh)
    \omega(\sqrt{2}\bh)\rd\bh=\frac{\la'!}{2^{\frac32}}\de_{\la',\bi},
\end{equation}
which can be computed using the orthogonality relation \eqref{eq:orth}. The proof is completed by substituting \eqref{eq:eta} into \eqref{eq:Aalk3}.
\end{proof}



\section{Proof of Lemma \ref{thm:int_Her}, Proposition \ref{thm:coe_D_psi} and \ref{thm:Maxwell}}
\label{app:VHS}
In this section, we provide the proofs to Lemma \ref{thm:int_Her}, Proposition \ref{thm:coe_D_psi} and \ref{thm:Maxwell}. 
{\renewcommand\proofname{Proof of Lemma \ref{thm:int_Her}}
\begin{proof}
From the recurrence relationship \eqref{eq:recur} of the Hermite polynomials, the one-dimensional standard Hermite polynomial can be expanded as 
\begin{equation}
    \label{eq:coef_Her}
    \begin{split}
        &H_{2n}(x)=\sum_{k=0}^n\frac{(2n-1)!!}{(2n-2k-1)!!}
        (-1)^kC_n^kx^{2n-2k}, \\
        &H_{2n-1}(x)=\sum_{k=0}^{n-1}\frac{(2n-1)!!}{(2n-2k-1)!!}
        (-1)^kC_{n-1}^kx^{2n-2k-1}.
    \end{split}
\end{equation}
Substituting \eqref{eq:coef_Her} into \eqref{eq:int_poly}, and let $\mC(\alpha_i,j_i), i = 1, 2,3$ be the coefficient of $x^{j_i}$ in $H_{\alpha_i}(x)$, we have
\begin{equation}
    \label{eq:int_poly_1}
     \mV(\kappa, \alpha, \mu) =   \sum_{\sk{\bj\pq\alpha \\ 2|(\alpha-\bj)}}\mC(\alpha_1, j_1)\mC(\alpha_2, j_2)\mC(\alpha_3, j_3)  \int_{\bbR^3}   \bv^{\kappa} \bv^{\bj}|\bv|^{\mu}\omega(\bv) \rd \bv.
\end{equation}
With the spherical coordinate transform $\bv=(r\chi_1,r\chi_2,r\chi_3)$, where $r\in \bbR^+$ and $\chi=(\chi_1, \chi_2, \chi_3)\in S^2$, it holds that 
\begin{equation}
    \label{eq:mV}
    \begin{split}
        \mathcal{V}(\kappa, \alpha, \mu)
        =&(2\pi)^{-\frac32}
        \sum_{\sk{\bj\pq\alpha \\ 2|(\alpha-\bj)}}
        \mC(\alpha,\bj) 
        \int_0^{\infty}r^{2+\mu+|\bj|+|\kappa|}\exp\left(-\frac{r^2}{2}\right)\rd r
        \int_{S^2} \chi_1^{j_1+\kappa_1}\chi_2^{j_2+\kappa_2}
        \chi_3^{j_3+\kappa_3}\rd\chi.
    \end{split}
\end{equation}
With Lemma \ref{thm:sp_int} and the properties of the Gamma function, \eqref{eq:mV} can be simplified as 
\begin{equation}
    \label{eq:mV_final}
    \begin{split}
        \mathcal{V}(\kappa, \alpha, \mu)
        =(2\pi)^{-\frac32}\sum_{\sk{\bj\pq\alpha \\ 2|(\alpha-\bj)}}
        \mC(\alpha, \bj)2^{\frac{1+\mu+|\bj|+|\kappa|}2} \Gamma\left(\frac{3+\mu+|\bj|+|\kappa|}2\right)\mS(\bj+\kappa).
    \end{split}
\end{equation}
This completes the proof. 
\end{proof}
}

{\renewcommand\proofname{Proof of Proposition \ref{thm:coe_D_psi}}
\begin{proof}
To prove \eqref{eq:lemma_D}, we first expand $H_{\alpha}(\bg')$ using \eqref{eq:coef_Her} as
\begin{equation}
    \label{eq:expan_Hg'}
    H_{\alpha}(\bg')=\sum_{\sk{\lambda\pq\alpha\\ 2|(\alpha-\lambda)}}\mC(\al_1,\la_1)\mC(\al_2,\la_2)\mC(\al_3,\la_3)(\bg')^{\lambda}. 
\end{equation}
Next, we use the binomial expansion for $(\bg')^{\lambda}$:
    \begin{equation}
        \label{eq:coe_D_1}
        (\bg')^{\lambda} =  \sum_{\sk{\kappa\pq\la }} C_{\lambda}^{\kappa}
         \left(\frac{1-e}{2}\right)^{|\kappa|}\left(\frac{1+e}{2}\right)
    ^{|\la|-|\kappa|} \bg^{\kappa} |\bg|^{|\lambda|-|\kappa|} \sigma^{\lambda - \kappa}. 
    \end{equation}
From Lemma \ref{thm:sp_int}, it can be observed that 
\begin{equation} 
\label{eq:coe_D_2}
\int_{S^2} \sigma^{\lambda - \kappa} \neq 0, \quad \text{ if and only if }  2|(\lambda - \kappa).
\end{equation}
Combining \eqref{eq:coe_D}, \eqref{eq:coe_D_1} and \eqref{eq:coe_D_2}, we can directly derive \eqref{eq:lemma_D}. For \eqref{eq:lemma_psi}, the proof is straightforward using the expansion of $H_{\alpha}(\bg)$ in \eqref{eq:coef_Her} and $\int_{S^2}1\rd\sigma=4\pi$.
\end{proof}
}

{\renewcommand\proofname{Proof of Proposition \ref{thm:Maxwell}}
\begin{proof}
    To prove this proposition, we start from \eqref{eq:coe_D} and \eqref{eq:coe_psi}. From \eqref{eq:coe_D_1} and \eqref{eq:coe_D_2}, we can derive that  $\bg^{\kappa}|\bg|^{|\lambda| - |\kappa|}$ is a polynomial of $\bg$ with degree $|\lambda|$. With the orthogonality of Hermite polynomials, it follows that 
         $$D(\bj, \kappa, 0)=0, \quad  \text{ if } |\bj|<|\kappa|. $$
    When $\mu = 0$, it is obvious that $\psi(\bj,\kappa,0)$ can be nonzero only when $\bj = \kappa$. Thus, one can see from \eqref{eq:cal_gamma} that 
     $$\gamma_{\kappa}^{\bj} =0, \quad  \text{ if } \varpi = 1~\text{and}~ |\bj|<|\kappa|. $$ 
     
     Finally, when $|\alpha| < |\lambda| + |\kappa|$ in $A\alk$, it can be observed that $|\kappa'|-|\bj|=|\kappa|+|\lambda|-|\alpha|>0$ in the summation \eqref{eq:theorem_A}. This completes the proof. 
\end{proof}
}

\section{Projection operator}
\label{app:project}
In this section, we present the theorem of the projection operator between different expansion centers. We refer the readers to \cite[Theorem 3.1]{ZhichengHu2019} for the related proof and details of this projection algorithm.
\begin{theorem}
    \label{thm:project}
    Suppose $f(\bv)$ 
    is expanded with two different expansion centers $[\ou\UP[1],\oT\UP[1]]$ and
    $[\ou\UP[2],\oT\UP[2]]$. From \eqref{eq:falpha}, we can compute the expansion coefficients for these two centers as
    \begin{equation}
        \label{eq:twoexpan}
        \begin{split}
            & f^{[\ou\UP[1],\oT\UP[1]]}_{\al}=\frac{1}{\al!}\int H_{\alpha}^{[\ou\UP[1],\oT\UP[1]]}(\bv)f(\bv)\rd\bv, \\
            & f^{[\ou\UP[2],\oT\UP[2]]}_{\al}=\frac{1}{\al!}\int H_{\alpha}^{[\ou\UP[2],\oT\UP[2]]}(\bv)f(\bv)\rd\bv. \\
        \end{split}
    \end{equation}
    Then we can obtain the second set of coefficients from the first set using the relationship
    \begin{equation}
        \label{eq:project_expan}
        f^{[\ou\UP[2],\oT\UP[2]]}_{\al}=\Big(\oT\UP[2]\Big)^{-\frac{|\al|}2}\sum_{l=0}^{|\al|}\phi_{\alpha}\UP[l], 
    \end{equation}
    where $\phi_{\alpha}\UP[l]$ is defined recursively as
    \begin{equation}
        \label{eq:def_phi}
        \phi_{\alpha}\UP[l]=\left\{
        \begin{array}{ll}
            \Big(\oT\UP[1]\Big)^{\frac{|\alpha|}2}f^{[\ou\UP[1],\oT\UP[1]]}_{\al}, & l=0,  \\
            \frac{1}{l}\sum_{d=1}^3\left[\Big(\ou\UP[2]-\ou\UP[1]\Big)\phi_{\al-e_d}\UP[l-1]+\frac12 \Big(\oT\UP[2]-\oT\UP[1]\Big)\phi_{\al-2e_d}\UP[l-1]\right], & 1\leqslant l \leqslant |\al|. 
        \end{array}
        \right.
    \end{equation}
    In \eqref{eq:def_phi}, terms with any negative index are 
    regarded as $0$.
\end{theorem}

\section{WENO reconstruction}
\label{app:WENO}
In this section, the WENO reconstruction for $\bdf$ is listed. The specific reconstruction coefficients are as follows: 
\begin{equation}
    \label{eq:WENO}
    \begin{split}
    & \bdf^{L,1}=\frac32\bdf^n_{j}-\frac12\bdf^n_{j-1}, \quad
    \bdf^{L,2}=\frac12\bdf^n_{j}+\frac12\bdf^n_{j+1}, \\
    & \bdf^{R,1}=\frac32\bdf^n_{j}-\frac12\bdf^n_{j+1}, \quad
    \bdf^{R,2}=\frac12\bdf^n_{j}+\frac12\bdf^n_{j-1}, \\
    & \om_{L,1}=\frac{\ga_1}{\Big[\eps+(\bdf^n_{j}-\bdf^n_{j-1})^2\Big]^2}, \quad
    \om_{L,2}=\frac{\ga_2}{\Big[\eps+(\bdf^n_{j+1}-\bdf^n_{j})^2\Big]^2}, \\
    & \om_{R,1}=\frac{\ga_1}{\Big[\eps+(\bdf^n_{j+1}-\bdf^n_{j})^2\Big]^2}, \quad
    \om_{R,2}=\frac{\ga_2}{\Big[\eps+(\bdf^n_{j}-\bdf^n_{j-1})^2\Big]^2}, \\
    &\bdf^{n,L}_{j+1/2}=\frac{\om_{L,1}\bdf^{L,1}+\om_{L,2}\bdf^{L,2}}
    {\om_{L,1}+\om_{L,2}}, \quad
    \bdf^{n,R}_{j-1/2}=\frac{\om_{R,1}\bdf^{R,1}+\om_{R,2}\bdf^{R,2}}
    {\om_{R,1}+\om_{R,2}}, \\
    & \eps=10^{-6},\quad \ga_1=\frac13, \quad \ga_2=\frac23,
    \end{split}
\end{equation}
where the square of $\bdf$ in \eqref{eq:WENO} indicates element-wise squaring, and the superscript $[\ou, \oT]$ on $\bdf$ is omitted.

\section{Nondimensionalization}
\label{app:nondim}
In this section, we provide the nondimensionalization to scale the variables as 
\begin{equation}
    \label{eq:nondim}
    \hat{\bx} = \frac{\bx}{x_0},  \qquad \hat{\bv} = \frac{\bv}{u_0}, \qquad     \hat{t}=\frac{t}{x_0/u_0}, \qquad \hat{m} = \frac{m}{m_0}, \qquad \hat{f} = \frac{f}{\rho_0 / (m_0 u_0^3)}, \qquad \hat{B}=\frac{B}{B_0},
\end{equation}
where $x_0, \rho_0$, and $m_0$ are the characteristic length, density and mass. Besides, $u_0$ is the character velocity defined as $u_0 = \sqrt{k_B \theta_0 / m_0}$ with $\theta_0$ the characteristic temperature. Besides, $B_0=\sqrt2u_0\pi d_{\rm ref}^2$ is adopted to rescale the HS collision kernel, where $d_{\rm ref}$ is the reference diameter.

The characteristic parameters are listed in Tab. \ref{table:character}, where $d_{\rm ref}$ is derived with \cite[Eq. (4.62)]{Bird}.


\begin{table}[!hptb]
\centering
\def\arraystretch{1.5}
{\footnotesize
\begin{tabular}{ll}
\hline
Characteristic parameters: & \\
Characteristic mass $m_0$ ($\times 10^{-26}$kg) & 6.63  \\
Characteristic length $x_0$ (m) & $10^{-3}$ \\
Characteristic velocity $u_0$ (m/s) & 238.377 \\
Characteristic temperature $\theta_0$ (K) & 273 \\
\hline 
Paramerters for HS model: & \\
Molecular mass: $m$ ($\times 10^{-26}$kg)          & 6.63           \\
Ref. viscosity: $\mu_{\rm ref}$ ($\times 10^{-5}$Pa s) & 2.117        \\
Viscosity index: $\varpi$              & 0.5 \\
Scattering parameter: $\alpha$ & 1 \\
Ref. diameter: $d_{\rm ref}$ ($\times 10^{-10}m$) & 3.63 \\
Ref. temperature: $T_{\rm ref}$ (K) & 273    \\
\hline
\end{tabular}
}
\caption{(Nondimensionalization in App. \ref{app:nondim}) Characteristic parameters in inhomogeneous tests.}\label{table:character}
\end{table}

\section{Properties of Hermite polynomials}
\label{app:Her}
For the Hermite polynomials \eqref{eq:Hermite}, several important properties are listed below:
\begin{property}
(Orthogonality)
\end{property}
\begin{equation}
    \label{eq:orth}
    \int_{\mR^3}H\aut(\bv)H\but(\bv)\omega_{\ou,\ot}(\bv)\rd\bv=\al_1!\al_2!\al_3!
    \delta_{\alpha,\beta}.
\end{equation}
\begin{property}
(Transitivity)
\end{property}
\begin{equation}
    \label{eq:tran}
    H\aut(\bv)=H_{\alpha}^{\bz,\zeta}\left(\sqrt{\frac{{\zeta}}{{\;\bT\;}}}(\bv-\ou)\right).
\end{equation}
\begin{property}
\label{property:recur}
(Recurrence)
\end{property}
\begin{equation}
    \label{eq:recur}
    \begin{split}
    &H^{\ou,\oT}_{\alpha+e_d}(\bv)=\frac{v_d-u_d}{\sqrt{\oT}}H^{\ou,\oT}
    _{\alpha}(\bv)-\al_dH^{\ou,\oT}_{\alpha-e_d}(\bv), \\
     & v_dH^{\ou,\oT}_{\alpha}=\sqrt{\oT}H^{\ou,\oT}_
    {\alpha+e_d}+
    u_dH^{\ou,\oT}_{\alpha}+\al_d\sqrt{\oT}H^{\ou,\oT}_{\alpha-e_d}.
    \end{split}
\end{equation}
\begin{property}(Differential of Hermite polynomial)
\label{property:deri}
\begin{equation}
    \label{eq:diff_Her}
    \frac{\pa}{\pa v_d}H\aut(\bv)=\frac{\al_d}{\sqrt{\oT}}H_{\alpha-e_d}^{\ou,\bT}
    (\bv).
\end{equation}
\end{property}
The property of transitivity can be directly derived from the definition \eqref{eq:Hermite}, and the proof of the other properties can be found in \cite{Abramowitz1964}.
\end{appendix}

\addcontentsline{toc}{section}{References}
\bibliographystyle{plain}
\bibliography{article}

\end{document}